\documentclass[11pt,a4paper]{article}

\usepackage{amsmath, amscd, amssymb, amsthm, MnSymbol}
\usepackage{hyperref}\usepackage{stmaryrd}
\usepackage{url}
\usepackage[all]{xy}
\usepackage{ulem}
\usepackage{tikz}
\usepackage{lineno}
\usepackage{arydshln}
\usepackage{circledsteps}

\allowdisplaybreaks

\newtheorem{theorem}{Theorem}[section]
\newtheorem*{theoremA}{Theorem A}

\newtheorem{proposition}[theorem]{Proposition}
\newtheorem{lemma}[theorem]{Lemma}

\newtheorem{question}[theorem]{Question}
\newtheorem*{claim}{Claim}
\newtheorem*{claimone}{Claim 1}
\newtheorem*{claimtwo}{Claim 2}
\newtheorem*{claimthree}{Claim 3}
\newtheorem*{claimfour}{Claim 4}
\newtheorem*{claimfive}{Claim 5}
\newtheorem*{claimsix}{Claim 6}
\newtheorem*{claimseven}{Claim 7}

\theoremstyle{definition} 
\newtheorem{note}[theorem]{Note}
\newtheorem*{notestar}{Note}
\newtheorem*{definitionstar}{Definition}
\newtheorem*{propositionstar}{Proposition}

\newtheorem{definition}[theorem]{Definition}
\newtheorem{remark}[theorem]{Remark}
\newtheorem{example}[theorem]{Example}
\newtheorem{notation}[theorem]{Notation}

\newcommand{\N}{{\mathbb N}}     
\newcommand{\Z}{{\mathbb Z}}     
\newcommand{\Q}{{\mathbb Q}}     
\newcommand{\RR}{{\mathbb R}}    
\newcommand{\Bool}{{\mathfrak{B}}}   %
\newcommand{\Latt}{{\mathfrak{L}}} 
\newcommand{\finind}{\textit{FinInd}}

\newcommand{\freez}{\textit{Frz}}

\newcommand{\suc}{\textit{Suc}}

\newcommand{\lcm}{\mathrm{lcm}}

\newcommand{\Suc}{\textit{Suc}}

\begin{document}
\setlength{\itemsep}{0pt}
\setcounter{page}{0}

\renewcommand\contentsname{\large \normalfont Title: {\bf
Congruence Preservation, Lattices and Recognizability}\\
Patrick Cegielski \& Serge Grigorieff \& Ir\`ene Guessarian}
\tableofcontents
\newpage

\title{Congruence Preservation, Lattices and Boolean Algebras}
\author{
{Patrick C\'egielski\textsuperscript{1}}
\and 
{Serge Grigorieff\textsuperscript{2} }
\and 
{Ir\`ene  Guessarian\textsuperscript{2,3}}
}
\maketitle

\footnotetext[1]{Emeritus at LACL, EA 4219, Universit\'e
Paris-Est Cr\'eteil, IUT S\'enart-Fontainebleau, France
\texttt{cegielski@u-pec.fr}.}
\footnotetext[2]{Retired from Universit\'e Paris 7 Denis
Diderot,
 France, \texttt {seg@irif.fr}, 
corresponding author.}
\footnotetext[3]{Emeritus at Sorbonne Universit\'e \texttt{ig@irif.fr}.}
\bibliographystyle{plain}

\maketitle

\begin{abstract} 
We study in general algebras Gr\"atzer's notion of
congruence preserving function, characterizing 
functions in terms of stability under inverse image of
particular Boolean algebras of subsets generated from any
subset of the algebra. 
Weakening Gr\"atzer's notion to only finite index
congruences, a similar result holds with lattices of sets.
Genereralizing the notion to that of stable preorder
preserving function, we extend these characterizations to
Boolean algebras and lattices generated from any
recognizable subset of the algebra.
Our starting point is a result with related flavor on the
additive algebra of natural integers which was obtained
some years ago.  
All these results can be visualized in the diagram of 
Table~\ref{table equivalences}.
We finally consider some simple particular conditions on
the algebra allowing to get a richer diagram.
\end{abstract}


\section{Introduction}
\label{s:intro}
%
The aim of this paper is to tackle the following question.
\begin{question}\label{question}
Is it possible to view a particular result about the particular
semiring  $\langle \N; +, \times \rangle$ (stated as Theorem
5.1 in \cite{cgg14ipl}) as an instance of a result about
general algebraic structures?
\end{question}

\subsection{Simple notions involved in the particular result
on $\langle \N; +, \times \rangle$}
\label{ss:motivation}
%
The notion of recognizability 
was originally introduced in the context of monoids of
words, cf. Sch\"utzenberger, 1956 \cite{Schutzenberger56}
p. 11, then in general algebras by Mezei \& Wright, 1967
\cite{MezeiWright67} (Definition 4.5, p. 22).
\begin{definition}\label{def recognizable}
A subset $X$ of an algebra $\+A$ is
$\+A$-\emph{recognizable} if $X=\varphi^{-1}(Z)$ for some
morphism $\varphi \colon \+A \to \+B$ into a finite algebra
and some subset $Z\subseteq \+B$. 
\end{definition}
Other used notions are about lattices.
\begin{definition}\label{def standardly bounded lattice}
1. A {\em lattice} $\Latt$ {\em of subsets} of a set $E$ is a
family of subsets of $E$ such that $L \cap M$ and
$L \cup M$ are in $\Latt$ whenever $L, M$ are in $\Latt$.

2. A lattice $\Latt$ of subsets of $E$ is {\em standardly
bounded} if $\emptyset$ and $E$ are in $\Latt$.
\end{definition}
\begin{definition}\label{def closed under inverse image}
A lattice $\Latt$ is {\em closed under inverse image by a
function} $g \colon E \to E$ if whenever it contains a set
$X$ it also contains the set
$g^{-1}(X) = \{ z \in E \mid g(z) \in X\}$.
In other words, $\Latt$ is {\em closed under} the function
$g^{-1}\colon \+P(E) \to \+P(E)$, where
$\+P(E)$ denotes the family of all subsets of $E$.
\end{definition}
\begin{notation}\label{not Latt suc (L)}
Given a set $L\subseteq\N$, we denote $\Latt_{\suc}(L)$
the smallest lattice of subsets of $\N$ which is closed under
the decrement function
$\suc^{-1}\colon \+P(\N) \to \+P(\N)$ and such that
$L \in \Latt$.
\end{notation}
We can now restate as Theorem~A the result of
our paper \cite{cgg14ipl} (Theorem 5.1) with two minor
modifications.

(1) First, we replace the notion of rational subset of $\N$
used in \cite{cgg14ipl} by that of recognizable subset of
$\N$, the two notions being equivalent on free monoids
(Kleene's theorem,  which we recall below in
\ref{ss rat rec}).

(2) We add an extra fifth condition
$(\Latt_{\suc}(L))^{\textit{rec}}_{\N}$ 
which is a straightforward avatar of the first condition
$(\textit{Latt}_{\suc})^{\textit{rec}}_{\N}$
and is similar to those dealing with finite sets and with
arithmetic progressions.
\begin{theoremA}
\label{thm origin}
Let $f \colon \N \to \N$. 
The following conditions are equivalent.
\\
$
\begin{array}{ll}
(\textit{Latt}_{\suc})^{\textit{rec}}_{\N}
&\textit{\begin{tabular}{l}
Every lattice $\Latt$ of recognizable subsets of $\N$ closed
under \\
decrement, i.e. under $\suc^{-1}$, is also closed under
$f^{-1}$
\end{tabular}}
\\
(\Latt_{\suc}(\forall L))^{\textit{rec}}_{\N}
&\textit{$\Latt_{\suc}(L)$ is closed under $f^{-1}$
for every recognizable subset $L \subseteq \N$}
\\
(\Latt_{\suc}(\forall L))^{\textit{fin}}_{\N}
&\textit{$\Latt_{\suc}(L)$ is closed under $f^{-1}$ 
for every finite subset $L \subseteq \N$}
\\
(\Latt_{\suc}(\forall L))^{\textit{arith}}_{\N}
&\left\{\textit{\begin{tabular}{l}
$\Latt_{\suc}(L)$ is closed under $f^{-1}$ for every  \\
arithmetic progression $L = q + r\N$ with $r > 0$
\end{tabular}}\right.
\\
(\textit{Arith})^{(abc)}
&\left\{\begin{array}{ll}
(a) &\textit{$|y-x|$ divides $|f(y) - f(x)|$ for all
$x, y \in \N$,}\\
(b) &\textit{$f(x) \geq x$ for all $x \in \N$} \\
(c) &\textit{$f$ is monotone non-decreasing}
\end{array}\right.
\end{array}
$
\end{theoremA}
\begin{notestar}
The reason for the a priori strange notation
$(\textit{Arith})^{(abc)}$ is because we also consider variant
of clause $(b)$, namely $(\textit{Arith})^{(ab^\flat)}$ and
$(\textit{Arith})^{(ab^\flat c)}$.
\end{notestar}
Let's just prove the equivalence between the first condition
--~which comes from \cite{cgg14ipl}~-- and the added fifth
condition $(\Latt_{\suc}(\forall L))^{\textit{rec}}_{\N}$.
\begin{proof}[Proof of 
$(\textit{Latt}_{\suc})^{\textit{rec}}_{\N} \Leftrightarrow
( \Latt_{\suc}(\forall L))^{\textit{rec}}_{\N}$]
${}$

For the left to right implication, observe that 
$\Latt_{\suc}(\forall L)$ consists of recognizable sets since
they are $(\cup, \cap)$-combinations of sets $\suc^{-n}(L)$,
with $n \in \N$, and recognizability is preserved under
$\suc^{-1}$, union and intersection.

For the right to left implication, let $\Latt$ be as in condition
$(\textit{Latt}_{\suc})^{\textit{rec}}_{\N}$ and $L \in \Latt$.
By condition $(\Latt_{\suc}(\forall L))^{\textit{rec}}_{\N}$, we
get $f^{-1}(L) \in \Latt_{\suc}(L)$, and inclusion
$\Latt_{\suc}(L) \subseteq \Latt$ insures that
$f^{-1}(L) \in \Latt$.
\end{proof}
\begin{remark}\label{rk surprising cong preserv N to N}
The successor function being such an elementary function,
it is tempting to imagine that the equivalent conditions of
Theorem~A entail that $f$ should itself be quite
simple\ldots but it turns out that this is far to be the case.

Indeed, functions satisfying clause $(a)$ of condition
$(\textit{Arith}^{(abc)})$ in Theorem~A are  characterized
in our paper \cite{cgg14}, and some of these functions
$\N \to \N$ are quite surprising, cf. Theorems 3.1 \& 3.13
\& Example 3.14 in \cite{cgg14} (as usual, $\lfloor z \rfloor$
is the integral part of a real $z$, and $e$ is Euler number,
and $\cosh$ and $\sinh$ are the hyperbolic cosine and sine
functions), for instance
\begin{equation*}
\begin{array}{l}
x\mapsto\left\{\begin{array}{ll}
1&\textit{if $x=0$}\\
 \lfloor e x! \rfloor&\textit{if $x\geq1$}
 \end{array}\right.
 \quad,\quad
x\mapsto\left\{\begin{array}{ll}
 \lfloor\cosh(1/2) 2^{x} x!\rfloor&\textit{if $x\in2\N$}\\
 \lfloor\sinh(1/2) 2^{x} x!\rfloor&\textit{if $x\in2\N+1$}
 \end{array}\right.
 \smallskip
 \\
 x \mapsto \lfloor e^{1/a} a^{x }x!\rfloor
 \qquad\textit{with $a \in \Z\setminus\{0,1\}$}
\end{array}
\end{equation*}
In particular, since $e = 2.718\ldots$, the first of these
functions is monotone non-decreasing and also satisfies
clause $(b)$ of condition $\textit{Arith}$, hence it satisfies
all conditions in Theorem~A.

In the ring $\langle  \Z; +, \times \rangle$ the analog of
condition $(\textit{Arith})$ is also be satisfied by very
intricate functions, for instance, cf. Theorem 26 in our
paper \cite{cgg14B}:
\begin{equation*}
f(n) = \left\{
\begin{array}{ll}
\sqrt{\dfrac{e}{\pi}}
\times \dfrac{\Gamma(1/2)}{2 \times 4^n\times n!}
\displaystyle\int_1^\infty e^{-t/2}(t^2-1)^n dt
&\quad\text{for $n\geq 0$}
\\
-f(|n|-1)&\quad\text{for $n<0$}
\end{array}\right.\,,
\end{equation*}
\end{remark}
%

%
\subsection{What is in the paper?}
\label{ss:roadmap}
%
\S\ref{s tools} presents tools used in this paper, most of
them being classical ones. 

Condition $(\textit{Arith}^{(abc)})$ of Theorem~A being
specific to arithmetic on $\N$, in order to get a version of
this theorem expressible in general algebras, we need to
modify it via some general algebraic notions.
It turns out that convenient notions are those of
{\em congruence preservation}, {\em stable-preorder
preservation} and {\em recognizabililty}.
The first two notions are recalled in \S\ref{ss gratzer}.

Recognizability is recalled in \S\ref{ss rec}

We explicit in \S\ref{ss reducing to arity 1} the freezification
of operations in a general algebra. This essentially reduces
the operations of the algebra to unary ones as concerns
what is studied in this paper.
This allows to introduce in \S\ref{ss syntactic} the syntactic
congruence and the syntactic stable preorder of a subset of
any general algebra, not only monoids contrarily to what
is always done, cf. \cite{almeida, pin95, pin22}.

\S\ref{s arith cong preserv} illustrates the notions of
\S\ref{s tools} in the algebras $\langle \N, Suc \rangle$,
$\langle \N, + \rangle$, and $\langle \N, +, \times \rangle$.

Folklore characterization of congruences and recognizability
are recalled in \S\ref{ss cong N} and \S\ref{ss rec N}.
Kleene's celebrated theorem about recognizability in free
monoids (which we use to go from our result in
\cite{cgg14ipl} to Theorem~A) is recalled in \S\ref{ss rat rec}.

In \S~\ref{ss:CP N} and \S\ref{ss:SPP N} condition
$(\textit{Arith}^{(abc)})$ of Theorem~A  --~which is
obviously specific to arithmetic on $\N$~--  is related to
the notions of congruence  and stable preorder
preservation in $\N$, (cf. Theorem~\ref{thm:CPsurN} and
Theorem~\ref{thm:miracle}).
These results give the key idea for extensions of
Theorem~A in a general algebra: replace condition
$(\textit{Arith}^{(abc)})$ --~which is specific to $\N$~-- by a
condition of preservation of congruences and stable
preorders  which makes sense in any algebra.
These results also illustrate the general algebraic fact that 
stable-preorder preservation may be strictly more
demanding than congruence preservation, 
cf. \S\ref{ss cong preorder preservation and monotonicity}.

In \S\ref{s lattice cong preserv} the notions of congruence
and stable preorder preservation in an algebra are related
to conditions on lattices and Boolean algebras much similar
to those of Theorem~A.
This is done via an explicit representation of $f^{-1}(L)$
--~in terms of the inverse images $\gamma^{-1}(L)$'s of the
freezifications~-- which leads to a lattice condition for
stable preorder preservation, and to a Boolean algebra
condition for congruence preservation
cf. \S\ref{ss f-1L preceq} and \S\ref{ss:preserv and latt BA}
and Theorem~\ref{th:LatticeGen}.
These conditions involve lattices and Boolean algebras
which are bounded when $L$ is recognizable (which is
always the case for Boolean algebras) and have to be
complete when $L$ is not recognizable.

In \S\ref{s improving} we consider diverse hypothesis on
the algebra to improve Theorem~\ref{th:LatticeGen},
replacing some or all implications by logical equivalences
in the figure of Table~\ref{table equivalences}.

\section{Algebraic tools used in this paper}\label{s tools}
%
\subsection{Gr\"atzer notion of congruence preservation in
general algebras}
\label{ss gratzer}
The first condition of Theorem A is related to the algebraic
notion of congruence preservation introduced by Gr\"atzer,
1962 \cite{gratzer62}  (under the name {\em substitution
property}). 

Let's recall this notion, together with a natural extension 
to {\em stable-preorder preservation} which, for
non-constant functions on $\N$, happens to be equivalent
to condition $(\textit{Arith})$ in Theorem~A, cf.
Theorem~\ref{thm:CPsurN} and
Theorem~\ref{thm:miracle}.
\begin{notation}\label{not algebra}
We denote $\+A = \langle A; \Xi \rangle$ the {\em algebra}
consisting of the nonempty carrier set $A$ equipped with a
set of operations $\Xi$, each $\xi \in \Xi$ being a mapping
$\xi \colon A^{ar(\xi)} \to A$, where $ar(\xi) \in \N$ is the
{\em arity} of $\xi$.
\end{notation}
\begin{definition}\label{def:fpreserveR}
Let  $\+A = \langle A; \Xi \rangle$ be an algebra.

1. A binary relation $\rho$ on $A$ is said to be
{\em compatible} with a function
$f \colon A^p \longrightarrow A$ if, for all elements
$x_1, \ldots, x_{p}, y_1, \ldots, y_{p}$ in $A$, we have
\begin{equation}\label{eq:preserveR}
(x_1 \rho\, y_1\ \wedge \cdots \wedge\ x_{p}
 \rho \, y_{p})\quad
  \Longrightarrow \quad f(x_1, \ldots, x_p) \; \rho \,
  f(y_1, \ldots, y_p).
\end{equation}  

2. A binary relation $\rho$ on $A$ is said to be
{\em $\+A$-stable} if it is compatible with each operation
$\xi \in \Xi$, i.e., \eqref{eq:preserveR} holds with $\xi$ in
place of $f$ for every $\xi \in \Xi$.
In particular, $\+A$-stable equivalence relations on $A$
are the {\em $\+A$-congruences}. 

3. A function $f \colon A^p \longrightarrow A$ is
{\em $\+A$-congruence preserving} if all
$\+A$-congruences are compatible with $f$, i.e.
\eqref{eq:preserveR} holds for every $\+A$-congruence
$\rho$.

4. Recall that a {\em preorder} is a reflexive and
transitive relation.
A function $f \colon A^p \longrightarrow A$ is
{\em $\+A$-stable-preorder preserving} if all
$\+A$-stable-preorders are compatible with $f$, i.e.
\eqref{eq:preserveR} holds for every
$\+A$-stable-preorder $\rho$.
\end{definition} 
As usual, when the algebra $\+A$ is clear from the context, 
$f$ is simply said to be {\em congruence preserving} 
(resp. {\em preorder preserving}).
\begin{remark}\label{rk composit preserves}
Clearly, any composition of functions in $\Xi$ preserves all
$\+A$-stable preorders of $\+A=\langle A;\Xi\rangle$.
\end{remark}
Since every congruence is a stable preorder, we get.
\begin{lemma}\label{l sp implies c}
If $f$ is $\+A$-stable preserving then it is also
$\+A$-congruence preserving.
\end{lemma}
\begin{remark}
The converse of Lemma~\ref{l sp implies c} is false: 
stable-preorder preservation is a strictly stronger condition
than congruence preservation since
Theorem~~\ref{thm:miracle} and
Proposition \ref{p:cp non monotone} below prove that there
exist functions which preserve
$\langle \N; + \rangle$-congruences but not all
$\langle \N; + \rangle$-stable-preorders.
\end{remark}
We shall also need the following classical notion.
\begin{definition}\label{def:finite index}
The {\em index of an $\+A$-congruence} is the number
(finite or infinite) of equivalence classes.
The {\em index of an $\+A$-stable-preorder} $\sqsubseteq$
is that of the {\em associated $\+A$-congruence}
$\{ (x, y) \mid x \sqsubseteq y \textit{ and }
 y \sqsubseteq x\}$.
\end{definition}

%
\subsection{Recognizable subsets of an algebra}
\label{ss rec}
%
Since congruences are kernels of morphisms, 
Definition~\ref{def recognizable} can be expressed in terms
of congruences.
\begin{lemma} \label{l:rec}\
A {\em subset} $L$ of an algebra $\+A$ is
$\+A$-{\em recognizable} if and only if it is saturated for
some finite index $\+A$-congruence (i.e. $L$ is a union of
equivalence classes).
\end{lemma}
We shall use the following straightforward result.
\begin{lemma}\label{l rec bool alg}
1. The family of recognizable subsets of $\+A$ is a Boolean
algebra.

2. If $L\subseteq \+B$ is $\+B$-recognizable and
$\varphi \colon \+A \to \+B$ is a morphism then 
$\varphi^{-1}(L)$ is $\+A$-recognizable.
\end{lemma}
We shall also use the following simple result.
\begin{lemma}\label{l rec initial segment}
Let $\sqsubseteq$ be a finite index $\+A$-stable preorder
on an algebra $\+A$.
Every $\sqsubseteq$-initial segment $L$ (i.e. whenever
$a \in L$ and $b \sqsubseteq a$ then $b \in L$) is
$\+A$-recognizable.
\end{lemma}
\begin{proof}
Observe that any initial segment of $\sqsubseteq$ is
saturated for the associated congruence which has finite
index. 
\end{proof}
%

%
\subsection{Freezifications}
\label{ss reducing to arity 1}
%
The notion of 1-freezification was not introduced for our
immediate problem but to characterize congruence
preservation of a function $f$ of arbitrary arity via
congruence preservation of suitably chosen unary functions
related to the function $f$.
It goes back to S\l{}omi\'nski, 1974, cf. Definition 1 in
\cite{Slominski74}.

This notion allows to show that, relative to congruences,
stable preorders and surjective homomorphisms, one can
replace the operations of an algebra by unary operations,
namely their freezifications.
\begin{definition}\label{not:freeze}
1. Let $f\colon A^{n}\to A$. 
If $n = 1$ the sole {\em 1-freezification} of $f$ is $f$ itself.
If $n \geq 2$ the {\em 1-freezifications} of
$f \colon A^{n} \to A$ are the unary functions
\begin{equation}\label{eq f i vec c}
f^{[i, \vec{c}]}(x) =
 f(c_1, \ldots, c_{i-1}, x, c_{i+1}, \ldots, c_{n})
\end{equation}
for $i \in \{1, \ldots, n \}$ and
$\vec{c} = (c_1, \ldots, c_{i-1}, c_{i+1}, \ldots, c_n)
 \in A^{n-1}$,
i.e., all but one argument of $f$ are frozen.
The set of 1-freezifications of $f$ is denoted $\freez(f)$.

2. If $\+F$ is a family of operations over a set $A$ with
arbitrary arities, the 1-freezification $\freez(\+F)$ of $\+F$
is the family of all 1-freezifications of operations in $\+F$.

3. A composition of a finite sequence of 1-freezifications is
called a {\em $*$-freezification}; in particular, the
composition of an empty sequence is the identity on $A$. 
The family of all $*$-freezifications of all operations in $\+F$
is denoted by $\freez^{*}(\+F)$.

4. To $\+A = \langle A; \Xi \rangle$ we associate two related
algebras $\freez(\+A)$ and $\freez^{*}(\+A)$ :
\begin{equation}\label{def:gen}
\freez(\+A) = \langle A; \freez(\Xi) \rangle
\qquad \textit{ and } \qquad
\freez^{*}(\+A) = \langle A; \freez^{*}(\Xi) \rangle.
\end{equation}
\end{definition}
\begin{figure}[h]
\center
\fbox{
$\begin{array}{lrcl}
\fbox{\textit{Over $\N$}}&\freez(\suc)&=&\{\suc\}\\
&{\freez^{*}(\suc), \freez(+), \freez^{*}(+)},&=&
\{\textit{translations }x\mapsto x+n \mid n\in\N\}\\
&{\freez(\times) = \freez^{*}(\times)}
&=&\{\textit{homotheties }x\mapsto nx \mid n\in\N\}\\
&\freez(\{+,\times\})
&=&\{\textit{translations and homotheties}\}\\
&\freez^{*}(\{+,\times\})
&=&\{\textit{affine maps }x\mapsto nx+p \mid n,p\in\N\}
\\
\cline{2-4}
\fbox{\textit{Over $\Sigma^*$}}&\freez(\cdot)
&=& \{x\mapsto ux,\ x\mapsto xv \mid u,v\in\Sigma^*\}
\\
&{\freez^{*}(\cdot)}
&=& \{x\mapsto uxv \mid u,v\in\Sigma^*\}
\end{array}$}
\caption{Examples of $\freez(\Xi)$, $\freez^{*}(\Xi)$}
\label{fig-freez}
\end{figure}
\begin{example}\label{ex freez}
Figure \ref{fig-freez} gives algebras $\freez(\+A)$ and
$\freez^{*}(\+A)$ for some algebras with carrier set $\N$
and for free monoids with the operation $\cdot$ of
concatenation.
\end{example}
\begin{remark}
In general, neither $\freez(\+A)$ nor $\freez^{*}(\+A)$
coincide with what are usually called ``affine'' or
``polynomial'' functions of the algebra (which should more
adequately be called ``term'' functions).

For instance, with $\langle \N; + \rangle$,
$\langle \N; \times \rangle$, and
$\langle \N; +, \times \rangle$ the usual ``affine'' or
``polynomial'' functions are respectively the usual affine
functions, monomial functions $x \mapsto nx^{p}$, and the
polynomial functions in $\N[x]$.

With $\langle \Sigma^*, \cdot \rangle$ the usual
``polynomial'' functions are the functions of the form
$x \mapsto u_{0}xu_{1}xu_{2} \ldots xu_{n}$, the $u_{i}$'s
being arbitrary words. If there are several occurrences of
$x$ in such a polynomial then it is neither in
$\freez(\langle \Sigma^*, \cdot \rangle)$ nor in
$\freez^*(\langle \Sigma^*, \cdot \rangle)$.
\end{remark}
The following result is inspired by the core of Theorem 1 in
S\l{}omi\'nski \cite{Slominski74} (trivially extending it from
equivalences to preorders). It shows that, relative to
congruences, stable preorders, surjective homomorphisms,
recognizability, one can replace the operations of the
algebra by unary operations, namely their freezifications.
\begin{lemma}\label{lem:de n a 1}
A preorder $\preceq$ on $A$ is compatible with
$f \colon A^n \longrightarrow A$ if and only if it is
compatible with all its 1-freezifications.
In particular, the three algebras $\+A$, $\freez(\+A)$, and
$\freez^*(\+A)$ have the same congruences and the same
stable preorders, hence the same functions preserving
congruences or stable preorders and the same recognizable
sets.
\end{lemma}
\begin{remark}
This possibility of reduction to arity 1 algebraic operations 
for key notions relative to algebras somewhat reminds of
the fact that in the framework of real analysis unary
continous functions and addition suffice to get all
continuous functions of any arity, more precisely 
$f(x_1,\ldots,x_{n}) 
= \sum_{p=1}^{p=2n+1} g_{p}(\sum_{q=1}^{q=n} h_{p,q}(x_{q}))$
(Andrei N. Kolmogorov solution of Hilbert's 13th problem,
1957, improving Vladimir Arnold solution).
\end{remark}
\begin{proof}[Proof of Lemma~\ref{lem:de n a 1}]
The left to right implication is a trivial use of the reflexivity
of $\preceq$. 
The right to left implication uses the transitivity of
$\preceq$. 
Indeed, if $\preceq$ is compatible with the 1-freezifications
\begin{equation*}
f^{[1,a_2,\ldots,a_{n}]},
f^{[2,b_1,a_{3}\ldots,a_{n}]},
f^{[3,b_1,b_{2},a_{4},\ldots,a_{n}]},
\ldots,
f^{[n-1,b_1,b_{2},\ldots,b_{n-2},a_{n}]},
f^{[n,b_1,b_{2},\ldots,b_{n-1}]}
\end{equation*}
and if $a_{1}\preceq b_{1}$,\ldots, $a_{n}\preceq b_{n}$ we
have
\begin{multline*}
f(a_1,a_2,\ldots,a_{n})
= f^{[1,a_2,\ldots,a_{n}]}(a_{1})
\preceq f^{[1,a_2,\ldots,a_{n}]}(b_{1})
= f(b_{1},a_2,\ldots,a_{n})
\\
= f^{[2,b_{1},a_{3},\ldots,a_{n}]}(a_{2})
\preceq f^{[2,b_{1},a_{3},\ldots,a_{n}]}(b_{2})
= f(b_{1},b_2,a_{3},\ldots,a_{n})
\\
= \ldots 
= f^{[n,b_{1},b_{2},\ldots,b_{n-1}]}(a_{n})
\preceq f^{[n,b_{1},b_{2},\ldots,b_{n-1}]}(b_{n})
= f(b_{1},b_2,\ldots,b_{n-1},b_{n})
\end{multline*}
hence $\preceq$ is compatible with $f$. 
\end{proof}
Finally, let's complete Lemma~\ref{l rec bool alg} as follows.
\begin{lemma}\label{l rec freez}
The family of recognizable subsets of $\+A$ is closed under
the $\gamma^{-1}(L)$'s for $\gamma \in \freez^{*}(+A)$.
\end{lemma}
\begin{proof}
If $L$ is saturated for a finite index congruence $\sim$ and
$x \in \gamma^{-1}(L)$, i.e. $\gamma(x) \in L$, the condition
$x\sim y$ entails
$\gamma(x) \sim \gamma(y)$ hence $\gamma(y) \in L$
since $L$ is saturated, i.e. $y \in \gamma^{-1}(L)$.
Thus, $\gamma^{-1}(L)$ is also saturated for $\sim$.
Therefore it is also recognizable.
\end{proof}

\subsection{Syntactic preorder and congruence of a subset
of an algebra}
\label{ss syntactic}
%
The notion of syntactic congruence first appeared in
Sch\"utzenber\-ger, 1956 \cite{Schutzenberger56} p.~11, in
the framework of the monoid of words with concatenation
over some alphabet (see also Mezei \& Wright, 1967
\cite{MezeiWright67}, S\l{}omi{\'n}ski \cite{Slominski74}).
It was extended to general monoids in an obvious way (see
Almeida \cite{almeida}, Pin, 2022 \cite{pin22}). 

The syntactic preorder on monoids appears as such in Pin,
1995 \cite{pin95}.

Using freezifications, we present these notions in the
framework of any algebra.
\begin{definition}\label{def syntactic cong pre}  
Let $\+A = \langle A; \Xi \rangle$ be an algebra. With any
$L \subseteq A$ are associated two relations $\preceq_L$
and $\sim_{L}$, called {\em syntactic preorder} and
{\em syntactic congruence} of $L$:
\begin{equation}\label{eq syntactic cong preord}
\begin{array}{rclrcl}
x\preceq_L y  &\text{ if and only if }&
\forall \gamma \in \freez^{*}(\+A)
& (\gamma(y) \in L& \Longrightarrow& \gamma(x) \in L)
 \\
x\sim_L y  &\text{ if and only if }&
\forall \gamma \in \freez^{*}(\+A)
&(\gamma(y) \in  L &\iff& \gamma(x) \in L)
\end{array}
\end{equation}
\end{definition}
\begin{lemma}\label{l preceqL stable}
1. Relations $\preceq_L$ and $\sim_L$ are an $\+A$-stable
preorder and an $\+A$-con\-gru\-ence. 

2. The congruence
$\{(x,y) \mid x \preceq_L y \textit{ and }y \preceq_L x\}$
is $\sim_L$.
\end{lemma}       
\begin{proof}
1. Apply Lemma~\ref{lem:de n a 1} to all operations in  $\Xi$.
Item 2 is obvious.
\end{proof}
\begin{example}
Let $\+N = \langle \N, \suc \rangle$.
The syntactic congruence of $L = a+k\N$ satisfies
\begin{align}\label{ex eq synt a+kN}
x \sim_L y &\Leftrightarrow \forall n \in \N\
\big(x+n \in (a+k\N) \Leftrightarrow y+n \in (a+k\N)\big)
\end{align}
Observe that $x \sim_L y$ implies $x \equiv y \bmod k$
since, letting $n$ be such that $x+n\in a+k\N$, we then
have $y+n\in a+k\N$ hence
$x-y \in (a+k\N)-(a+k\N) = k\Z$.

We show that $\sim_{L}$ coincides with
$\equiv_{\max(0,a-k+1),k}$ (defined by \eqref{equation:ak}
in Lemma~\ref{l folklore cong}).
We reduce to prove that they coincide on pairs $(x, y)$ such
that $x < y$.
\\
\textbullet\
Assume $x \equiv_{\max(0,a-k+1),k} y$. We prove that
$x \sim_L y$.   
Since $x < y$, we have $\max(0,a-k+1) \leq x$ and
$x \equiv y \bmod k$.
Implication $x+n \in (a+k\N) \Rightarrow y+n \in (a+k\N)$
is clear since $x \equiv y \bmod k$ and $x<y$ .
In case $y+n \in (a+k\N)$, since $x \equiv y \bmod k$
we have $x+n \in a+k\Z$, say $x+n=a+kz$ with $z \in \Z$.
Now, $x+n \geq \max(0, a-k+1+n)$ hence
$a+kz \geq a-k+1$, i.e. $k(z+1)\geq 1$ hence $k \in \N$
and $x+n \in a+k\N$.
This proves $x \sim_L y$.   
\\
\textbullet\
Assume $x \not\equiv_{\max(0, a-k+1),k} y$. We prove that
$x \not\sim_L y$.  
Since $x < y$ there are two possible reasons for
$x \not\equiv_{\max(0, a-k+1),k} y$.
\\
{\it Case $x \not\equiv y \bmod k$.} Then $x \not\sim_L y$
is clear.
\\
{\it Case $x \equiv y \bmod k$ and $x < \max(0, a-k+1)$.}
Let $p \geq 1$ be such that $y = x+kp$.
Since $x \geq 0$ we must have $a-k+1 \geq 1$, i.e.
$a \geq k$, and also $x < a-k+1$, i.e. $a-x \geq k$.
Let $n = a-x-k \in \N$. We have $x+n = a-k\notin a+k\N$.
Now, $y+n=x+n+pk=a+k(p-1)\in a+k\N$.
Thus, $x \not\sim_L y$.
This concludes the proof that $\sim_{L}$ and
$\equiv_{\max(0,a-k+1),k}$ coincide.
\end{example}
\begin{proposition}\label{p:syntactic largest}
Let $L$ be a subset of an algebra $\+A$.

1. $L$ is $\sim_L$-saturated and is an $\preceq_L$-initial
segment. 

2. If $L$ is saturated for a congruence $\equiv$ of $\+A$
then $\equiv$ refines the syntactic congruence $\sim_L $ of
$L$, i.e., $x \equiv y$ implies $x \sim_L  y$.
In particular, $L$ is recognizable if and only if $\sim_L$ has
finite index and if and only if $\preceq_L$ has finite index.

3. If $L$ is an initial segment of a stable preorder $\preceq$
of $\+A$ then $\preceq$ refines the syntactic preorder
$\preceq_L $ of $L$, i.e., $x \preceq y$ implies
$x \preceq_L  y$.
\end{proposition}
\begin{proof}
1. Let $\gamma$ be the identity in
\eqref{eq syntactic cong preord}.

2. The usual proof in textbooks for the case of monoids
(e.g. Almeida \cite{almeida} Sections 0.2, 3.1, 
or Pin \cite{pin22}, chap. IV, \S4.1, 5.3)  
extends to any algebra using freezifications.
Suppose $x \equiv y$. 
Since $\equiv$ is an $\+A$-congruence, for all
$\gamma \in \freez^{*}(\+A)$ we have
$\gamma(x) \equiv \gamma(y)$.
Since $L$ is $\equiv$-saturated we have
$\gamma(y) \in L \Leftrightarrow \gamma(x) \in L$. 
Being true for all $\gamma \in \freez^{*}(\+A)$, this means
$x \sim_{L} y$. 
3. Similar, initial segments replacing saturated sets.
\end{proof}
Let's state a straightforward consequence of item 1 of
Proposition~\ref{p:syntactic largest}.
\begin{lemma}\label{l preceqL}
$L$ is recognizable if and only if $\sim_L$ has finite index.
\end{lemma}       
%

\section{Case of the algebra $\langle \N, \suc \rangle$ or
$\langle \N, + \rangle$ or $\langle \N, +, \times \rangle$}
\label{s arith cong preserv}
%
%
\subsection{Folklore characterization of congruences on
$\+N$}
\label{ss cong N}
%
To show how Theorem~A is an instance of our general
results in \S\ref{s lattice cong preserv}, we need the classical
characterization of congruences on $\N$.

First, a straightforward observation.
\begin{lemma}\label{l folklore suc + +times idem}
The three algebras 
$\langle \N; \Suc \rangle$,
$\langle \N; + \rangle$
and $\langle \N; +, \times \rangle$
have the same stable preorders and the same congruences.
A fortiori, they admit the same stable-preorder-preserving
functions and the same  congruence-preserving functions.
\end{lemma}
\begin{proof}
Stability of a relation $\rho \subseteq \N\times\N$ under
successor implies stability under addition: to deduce
$(x+z) \rho (y+z)$ from $x \rho y$ merely apply $z$
times stability under successor.
Stability of a relation $\rho \subseteq \N\times\N$ under
addition implies stability under multiplication: straigtforward
induction on $t$ since from $x \rho y$ and $(xt) \rho (yt)$
we deduce $(xt+x) \rho (yt+y)$, that is
$x(t+1) \rho y(t+1)$.
\end{proof}
\begin{remark}\label{rk Suc + times preserv idem}
A priori, $\+A$-congruences and $\+A$-stable-preorders
strongly depend upon the signature of $\+A$. 
However, Lemma~\ref{l folklore suc + +times idem}
witnesses that this is not always the case.
\end{remark}
\begin{definition}
A {\em congruence} $\sim$ is {\em left cancellable} if
$\xi(a, x) \sim \xi(a, y)$ implies $x \sim y$ for each binary
operation $\xi$ of the algebra. {\em Right cancellability}
and {\em cancellability} are similarly defined.
\end{definition}
Congruences for the usual structures on $\N$ are simply
characterized, cf. Preston, 1963, \cite{HowiePreston}
(Lemma~2.15 \& Remark pp.~28-29),
Grillet, 1995 \cite{Grillet} (Lemma~5.6 p.~19).
\begin{lemma}\label{l folklore cong}
1. The congruences of the algebras 
$\langle \N; \Suc \rangle$,
$\langle \N; + \rangle$,
and $\langle \N; +, \times \rangle$
are the identity relation and the $\equiv_{a, k}$'s for
$a, k \in \N$, $k \geq 1$, where
\begin{equation}\label{equation:ak}
x \equiv_{a, k} y \mbox{  if and only if  }
\begin{cases}
\mbox{either } x = y\\
\mbox{or }  a \leq x 
\mbox{ and } a \leq y
\mbox{ and }  x \equiv y \pmod k
\end{cases}
\end{equation}
Congruence $\equiv_{a, k}$ has finite index $a+k$, i.e.
has $a + k$ equivalence classes.
It is cancellable if and only if $a = 0$.

2. Up to isomorphism, the quotient algebra 
$\langle\N; \suc, +, \times\rangle / \equiv_{a, k}$
can be seen as the algebra
$\langle\{0, \ldots, a+k-1\};
 \Circled{suc}, \oplus, \otimes \rangle$
where 
\begin{align*}
\Circled{suc}(x) &= \texttt{IF } \suc(x) < a+k \texttt{ THEN }
 \suc(x) \texttt{ ELSE } a
\\
x\oplus y &= \texttt{IF } x+y < a+k \texttt{ THEN } x+y
\texttt{ ELSE } a + ((x+y-a) \bmod k)
\\
x\otimes y &= \texttt{IF } x \times y < a+k \texttt{ THEN }
x \times y \texttt{ ELSE } a + ((x \times y -a) \bmod k)
\end{align*}
\end{lemma}

%
\subsection{Folklore characterization of recognizable
subsets of $\+N$}
\label{ss rec N}
%
First, observe the following immediate consequence of 
Lemma~\ref{l folklore suc + +times idem}.
\begin{lemma}\label{l folklore suc + +times idem for rec}
The three algebras 
$\langle \N; \Suc \rangle$,
$\langle \N; + \rangle$,
and $\langle \N; +, \times \rangle$
have the same recognizable subsets.
\end{lemma}
\begin{proof}
By Lemma~\ref{l folklore suc + +times idem} they have the
same congruences hence the same subsets saturated for
some finite index congruence.
\end{proof}

In case the monoid is $\langle \N; + \rangle$, the
characterization of recognizable sets dates back to
Myhill, 1957 \cite{Myhill57}, cf. also
Eilenberg~\cite{Eilenberg74}, p. 101 Proposition 1.1.
The published proofs are (as far as we know) always
presented in terms of automata, we present it in terms of
congruences.
\begin{lemma}\label{l rec N}
Let $L \subseteq \N$. 
The following conditions are equivalent.
\begin{enumerate} \setlength\itemsep{0em}
\item
$L$ is a recognizable subset of the monoid
$\langle \N, + \rangle$.
\item
$L$ is the union of a finite set with a finite (possibly empty)
family of arithmetic progressions.
\item
$L$ is finite or $L = A \cup (a+\N)$ with $A$ a strict subset
of $\{0, \ldots, a-1\}$ possibly empty or
$L = A \cup (\Delta+k\N)$ with $k \geq 2$ and $A$ a
possibly empty subset of $\{0, \ldots, k-1\}$ and $\Delta$ a
non-empty strict subset of $\{0, \ldots, k-1\}$.
\end{enumerate}
\end{lemma}
\begin{proof}
$(1)\Rightarrow(3)$.
Any recognizable subset $L$ of $\N$ is
$\sim_{a,k}$-saturated for some $a \in \N$ and $k \geq 1$,
hence is of the form $L = A \cup ((a+\Delta)+k\N)$ for some
possibly empty sets $A \subseteq \{0, \ldots, a-1\}$ and
$\Delta \subseteq \{0, \ldots, k-1\}$.
If $k = 1$ we get $L = A$ or $L = A \cup (a+\N)$.
If $k \geq 2$ observe that the equality
$\Delta = \{0, \ldots, k-1\}$ implies
$(a+\Delta)+k\N = a+\N$.

$(3)\Rightarrow(2)$ is trivial.

$(2)\Rightarrow(1)$. An arithmetic progression $a + k\N$ is
the $\sim_{a, k}$ equivalence class of $a$, hence is
recognizable.
A singleton set $\{b\}$ is $\sim_{a, 1}$-saturated for any
$a > b$, hence is recognizable.
To conclude use closure under union of recognizability.
\end{proof}

\subsection{Kleene's theorem: recognizable and rational
sets in free monoids}
\label{ss rat rec}
%
First, recall a classical notion from theoretical computer
science, cf. for instance
Eilenberg \& Sch\"utzenber\-ger, 1969 
\cite{EilenbergSchutz69}, 
or Eilenberg, 1974 \cite{Eilenberg74} p. 159 sq.).

\begin{definitionstar}\label{def reg}
Let $\langle M, \cdot \rangle$ be a monoid with neuter
element $\varepsilon$. The family of {\em rational} (also
called {\em regular}) subsets of $M$ is the smallest family
$\+R$ of subsets of $M$ such that
(1) the empty set and every singleton set is in $\+R$,
(2) if $P$ and $Q$ are in $\+R$ then so are the union
$P \cup Q$, the product
$PQ = \{ u\cdot v \mid u \in P, v \in Q\}$, and the infinite
union of powers (also called Kleene's star)
$P^* = \bigcup_{n \, \in \, \N}P^{n}$, where
$P^{0}=\{\varepsilon\}$ and $P^{n+1}=P^{n}P$.
\end{definitionstar}
A celebrated theorem due to Kleene, 1956 \cite{Kleene56}, 
relates rationality and recognizability in a free monoid, i.e.
(up to isomorphism) an algebra of words over some
alphabet endowed with the concatenation operation,
cf. Eilenberg \cite{Eilenberg74} (p. 175 Theorem 5.1)
or Almeida~\cite{almeida} or Pin~\cite{pin22} (p. 70).
\begin{propositionstar}[Kleene's theorem] \label{p Kleene}
In free monoids, the notion of rational subset coincides with
that of recognizable subset.
In particular, this is true in the monoid
$\langle \N, + \rangle$ 
(integers being seen as words on a singleton alphabet).
\end{propositionstar}

\subsection{$\+N$-congruence preservation and arithmetic}
\label{ss:CP N}
%
Congruence preservation on $\langle \N, Suc \rangle$,
$\langle \N, + \rangle$, $\langle \N, +, \times \rangle$
can be simply characterized as follows with a weakening
$(\textit{Arith}^{(ab^\flat)})$ of condition
$(\textit{Arith}^{(abc)})$ of Theorem~A in
\S\ref{ss:motivation} which gives up the monotony of $f$
(clause $(c)$) and weakens clause $(b)$ by allowing $f$ to
be constant.
\begin{theorem}\label{thm:CPsurN}
For any $f \colon \N \longrightarrow \N$, the following
conditions are equivalent
\begin{equation*}
\begin{array}{ll}
(\textit{CongPres})_{\N}&
\textit{$f \colon \N \to \N$ is $\+N$-congruence
 preserving.}
\smallskip\\
(\textit{Arith}^{(ab^\flat)})&
\left\{\begin{array}{ll}
(a) &
\textit{$|y-x|$ divides $|f(y)-f(x)|$ for all $x, y \in \N$}\\
(b)^{\flat} &\textit{either $f$ is constant or
 $f(x) \geq x$ for all $x$}
\end{array}\right.
\end{array}
\end{equation*}
\end{theorem}
\begin{proof}
$(\textit{CongPres})_{\N} \Rightarrow
 (\textit{Arith}^{(ab^\flat)})$.
Suppose $f$ is congruence preserving.
Let $x < y$ and consider the congruence modulo $y-x$.
As $y \equiv x \bmod y-x$ we have
$f(y) \equiv f(x) \bmod y-x$, hence $y-x$ divides $f(y)-f(x)$.
This proves clause $(a)$ of condition
$(\textit{Arith}^{(ab^\flat)})$.
 
To prove clause $(b)^{\flat}$ of condition
$(\textit{Arith}^{(ab^\flat)})$, we assume that condition
$f(x) \geq x$ is not satisfied and show that $f$ is constant.
Let $a$ be least such that $f(a) < a$.
Consider the congruence $\equiv_{a, 1}$.
We have $a \equiv_{a, 1} y$ for all $y \geq a$, hence
$f(a) \equiv_{a, 1} f(y)$. 
As $f(a) < a$, this implies $f(y) < a$ and $f(a) = f(y)$.
Thus, $f(a) = f(y)$ for all $y \geq a$.
Let $z < a$. By clause (a) (already proved) we know that $p$
divides $f(z) - f(z+p)$ for all $p$.
Now, $f(z + p) = f(a)$ when $z + p \geq a$.
Thus, $f(z) - f(a)$ is divisible by all $p \geq a-z$.
This shows that $f(z) = f(a)$.
Summing up, we have proved that $f$ is constant with value
$f(a)$.

$(\textit{Arith}^{(ab^\flat)}) \Rightarrow
 (\textit{CongPres})_{\N}$.
Constant functions are trivially congruence preserving.
We thus assume that $f$ is not constant, hence $f$ satisfies
clause $(a)$ of condition $(\textit{Arith}^{\flat})$ and
$f(x) \geq x$ for all $x$.
Consider a congruence $\sim$ and suppose $x \sim y$. 
If $\sim$ is the identity relation then $x = y$ hence
$f(x) = f(y)$ and $f(x) \sim f(y)$.
Else, by Lemma~\ref{l folklore cong} the congruence $\sim$
is some $\equiv_{a, k}$ with $a \in \N$ and $k \geq 1$.
In case $y < a$ then condition $x \equiv_{a, k} y$ implies
$x = y$, hence $f(x) = f(y)$ and $f(x) \sim f(y)$.
In case $y \geq a$ then condition $x \equiv_{a, k} y$ implies
$x \geq a$ and $x \equiv y \bmod k$, hence $k$ divides
$y-x$.
Since clause $(a)$ of condition $(\textit{Arith}^{(ab^\flat)})$
ensures that $y-x$ divides $f(y)-f(x)$, we see that $k$ also
divides $f(y)-f(x)$.
Also, our hypothesis yields $f(x) \geq x$ and $f(y) \geq y$.
As $x, y \geq a$, we get $f(x), f(y) \geq a$. Since $k$
divides $f(y)-f(x)$, we conclude that $f(x) \equiv_{a, k} f(y)$,
whence $(\textit{CongPres})_{\N}$.
\end{proof}
\begin{remark}\label{exCPnonsurlineaire}
Clause $(a)$ of condition $(\textit{Arith}^{(ab^\flat)})$ does
not imply clause $(b)^{\flat}$, hence cannot be removed in
Theorem~\ref{thm:CPsurN}.
For instance, let $F \colon \N \to \N$ be the polynomial
$F(x) = x(x-1)\ldots(x-k)$ with $k\geq1$.
Since $y-x$ divides $y^{n}-x^{n}$ for any $n\geq1$,
developping $F(x)$ in a sum of monomials, we see that it
also divides $F(y)-F(x)$.
However, $F$ is neither constant nor overlinear as
$F(i) = 0 < i$ for $i = 1, \ldots, k$.
\end{remark}
%

%
\subsection{$\+N$-stable-preorder preservation and
arithmetic}
\label{ss:SPP N}
%
We now look at what monotonicity adds to congruence
preservation in $\langle \N; Suc \rangle$,
$\langle \N; + \rangle$, and $\langle \N; +, \times \rangle$,
namely, it insures stable-preorder preservation.
First, a simple observation.
\begin{lemma}\label{l:cp and spp on N}
Let $\preceq$ be a $\langle \N; + \rangle$-preorder and, for
$a \in \N$, let $M^+_a = \{ x \mid a \preceq a+x\}$ and
$M^-_a = \{ x \mid a+x \preceq a\}$.

\smallskip

1) $M^+_a$ and $M^-_a$ are submonoids of
$\langle \N; + \rangle$.

\smallskip

2) If $a \leq b$ then $M^+_a \subseteq M^+_b$ and
$M^-_a \subseteq M^-_b$.
\end{lemma}
\begin{proof}
1) Clearly, $0 \in M^+_a$. 
Suppose $a \preceq a+x$ and $a \preceq a+y$.
By stability, the first inequality yields $a+y \preceq a+x+y$
and, by transitivity, the second inequality gives
$a \preceq a+x+y$.
Thus, $M^+_a$ is a submonoid. Idem with $M^-_a$.

2) Since $b-a \in \N$, if $a \preceq a+x$ then by stability
$a+(b-a) \preceq a+x+(b-a)$, i.e., $b \preceq b+x$.
Thus, $M^+_a \subseteq M^+_b$.
Similarly, $M^-_a \subseteq M^-_b$. 
\end{proof}
\begin{theorem}\label{thm:miracle}
Let $\+N$ be any of the algebras 
$\langle \N; \Suc \rangle$ or
$\langle \N; + \rangle$ or
$\langle \N; +, \times \rangle$.
For any $f \colon \N \to \N$, the following conditions are
equivalent.
\begin{equation*}
\begin{array}{ll}
(\textit{PreordPres})_{\N}& \text{$f$ is
$\+N$-stable-preorder preserving}
\\
(\textit{CongPres}^{\uparrow})_{\N}& 
\text{$f$ is monotone non-decreasing and
$\+N$-congruence preserving}
\\
(\textit{Arith}^{(ab^\flat c)})&
\left\{\begin{array}{ll}
(a) &\textit{for all $x, y \in \N$, $y-x$ divides $f(y)-f(x)$}\\
(b)^{\flat} &\textit{$f$ is either constant or $f(x) \geq x$ for
all $x \in \N$}\\
(c) &\textit{$f$ is monotone non-decreasing.}
\end{array}\right.
\end{array}
\end{equation*}
\end{theorem}
\begin{proof}
$(\textit{PreordPres})_{\N} \Rightarrow
 (\textit{CongPres}^{\uparrow})_{\N}$.
If $f$ is stable preorder preserving then it is a fortiori
congruence preserving. 
As the usual order $\leq$ is $\+N$-stable, it is preserved by
$f$, hence $f$ is monotone non-decreasing.

$(\textit{CongPres}^{\uparrow})_{\N} \Leftrightarrow
 (\textit{Arith}^{(ab^\flat c)})$.
This is Theorem~\ref{thm:CPsurN} for monotone
non-decreasing $f$'s.

$(\textit{Arith}^{(ab^\flat c)}) \Rightarrow
 (\textit{PreordPres})_{\N}$.
Assume $(\textit{Arith}^{(ab^\flat c)})$.
We prove that $f$ preserves every stable preorder
$\preceq$.
The case where $f$ is constant is trivial. 
We now suppose $f$ is not constant.
Suppose $a \preceq b$, we argue by distinguishing two
cases.

{\it Case $a \leq b$.}
Then $b-a \in M^+_a = \{x \mid a\preceq a+x\}$.
By clause $(a)$ of condition $(\textit{Arith}^{(ab^\flat c)})$
we know that $b-a$ divides $f(b)-f(a)$. 
Now, $f(b)-f(a) \geq 0$ since $f$ is monotone
non-decreasing. As a consequence, $f(b)-f(a) = n(b-a)$ for
some $n \geq 0$. 
As $b-a \in M^+_a$ and $M^+_a$ is a monoid, $f(b)-f(a)$
is also in $M^+_a$.
Since $f(x) \geq x$ for all $x$, we have $f(a) \geq a$ hence,
by Lemma~\ref{l:cp and spp on N},
$M^+_a \subseteq M^+_{f(a)}$.
Thus, $f(b)-f(a)$ is in $M^+_{f(a)}$, implying
$f(a) \preceq f(b)$.

{\it Case $b \leq a$.}
Similar proof, observe that
$a-b \in M^-_b = \{ x \mid b+x \preceq b\}$.
\end{proof}
%

\subsection{Preservation of congruences may not extend to
stable preorders}
\label{ss cong preorder preservation and monotonicity}
%
Condition $(\textit{Arith}^{(abc)})$ of Theorem~A in
\S\ref{ss:motivation} and condition
$(\textit{CongPres}^{\uparrow})_{\N}$ of
Theorem~\ref{thm:miracle} ask for $f$ to be monotone
non-decreasing.
We show that this requirement cannot be removed: there
does exist congruence preserving functions over $\N$ which
fail monotonicity, cf. Proposition~\ref{p:cp non monotone}
below.

First, recall a result, cf. Theorem~3.15 in our paper
\cite{cgg14}.
\begin{proposition}[\cite{cgg14}]\label{p:g idr f}
For every $f \colon \N \to \N$ there exists a function
$g \colon \N \to \Z$ such that, letting $\lcm(x)$ be the least
common multiple of $1, \ldots, x$,
\begin{enumerate}
\item[(i)] 
$x-y$ divides $g(x)-g(y)$ for all $x, y \in \N$,
\item[(ii)] 
$f(x) - 2^x \lcm(x) \leq g(x)\leq f(x)$.
\end{enumerate}
Such a $g$ is
$g(x) = \sum_{k \in \N} \widetilde{a_k}\dbinom{x}{k}$,
where $f(x) = \sum_{k \in \N} a_k \dbinom{x}{k}$
is the Newton representation of $f$ (where the $a_{k}$'s are
in $\Z$) and $\widetilde{a_k} = \lcm(k)q_k$, where
$a_k = \lcm(k)q_k+b_k$ with $0 \leq b_k < \lcm(k)$.
\end{proposition}
\begin{proposition}\label{p:cp non monotone}
There exists an $\langle \N; +\rangle$-congruence
preserving function $g \colon \N \to \N$ which
is non-monotonous, hence is not
$\langle \N; +\rangle$-stable preorder preserving.
\end{proposition}
\begin{proof}
Consider the function $f$ such that
\begin{align}\notag
f(x) &= \sum \{2^{y+3}\lcm(y+2) \mid y\leq x, \ x\textit{ even}\}
          - \sum \{2^{z+1}\lcm(z) \mid z\leq x, \ z\textit{ odd}\}
\end{align}
We have $f(0) = 8\lcm(2) = 16$ and, for $p \in \N$,
\begin{align}\label{eq f(2p+2)-f(2p+1)}
f(2p+2) - f(2p+1) &= 2^{2p+5}\lcm(2p+4)
\\\label{eq f geq even}
f(2p)&\geq 2^{2p+3}\lcm(2p+2) > 0
\\\label{eq f(2p+1)-f(2p)}
f(2p+1)-f(2p) &= - 2^{2p+2}\lcm(2p+1)
\\\notag
f(2p+1)&\geq 2^{2p+3}\lcm(2p+2) - 2^{2p+2}\lcm(2p+1)
 \\\label{eq f geq odd}
             &\geq 2^{2p+2}\lcm(2p+1) > 0
 \\\label{eq f geq}
 \text{\eqref{eq f geq even} and \eqref{eq f geq odd} show
  that }
 & \text{$f$ maps $\N$ into $\N$ and
 $f(x) \geq 2^{x+1}\lcm(x)$}
\end{align}
Let $g \colon \N\to\Z$ be the function  associated to $f$
by Proposition~\ref{p:g idr f}.
We first prove that $g$ also maps $\N$ into $\N$. 
Clause (ii) in Proposition \ref{p:g idr f} and inequality
\eqref{eq f geq} insure that
$g(x) \geq f(x) - 2^x\lcm(x) > 2^{x} \, \lcm(x)$.
This proves that $g(x) \in \N$ and moreover $g(x) > x$ for
all $x$.

We next prove that $g$ is non monotone.
Using clause (ii) and \eqref{eq f geq}, we get
\begin{equation}\label{eqn:g idr f}
f(x+1) - f(x) - 2^{x+1} \lcm(x+1)  \leq g(x+1)-g(x) 
                                                    \leq f(x+1) - f(x) + 2^{x} \lcm(x)
\end{equation}
Using \eqref{eqn:g idr f} and \eqref{eq f(2p+1)-f(2p)} 
and \eqref{eq f(2p+2)-f(2p+1)}, we get
\begin{align}\notag
g(2p+1) - g(2p)& \leq f(2p+1)-f(2p) +2^{2p}\lcm(2p)
\\\label{eq g(2p+1)-g(2p)}
&= - 2^{2p+2}\lcm(2p+1) + 2^{2p}\lcm(2p) < 0
\\\notag
g(2p+2) - g(2p+1)& \geq f(2p+2) - f(2p+1)
 - 2^{2p+2}\lcm(2p+2) 
\\\label{eq g(2p+2)-g(2p+1)}
&=  2^{2p+5}\lcm(2p+4) - 2^{2p+2}\lcm(2p+2) > 0
\end{align}

Inequalities \eqref{eq g(2p+1)-g(2p)} and
\eqref{eq g(2p+2)-g(2p+1)} show that $g$ is a
non-monotone zigzag function.

Finally, as $g(x) > x$ for all $x$ and $x-y$ divides $g(x)-g(y)$
for all $x, y \in \N$ (by clause (i) of Proposition \ref{p:g idr f}), 
both clauses $(a)$ and $(b)^{\flat}$ of condition
$(\textit{Arith}^{(ab^\flat)})$ of Theorem \ref{thm:CPsurN}
hold, hence $g$ is congruence preserving.
\end{proof}

%
\subsection{A strengthening of Theorem~\ref{thm:miracle}}
%
Theorem~A in \S\ref{ss:motivation} relates conditions
involving particular lattices of subsets of $\N$ closed under
$f^{-1}$ with an arithmetic condition
$(\textit{Arith}^{(abc)})$ on $f$.
On its side Theorem~\ref{thm:miracle} relates stable
preorder preservation of $f$ with an arithmetic condition
$(\textit{Arith}^{(ab^\flat c)})$ which is equivalent to the
disjunction of condition $(\textit{Arith})^{(abc)}$ of
Theorem~A and the condition that $f$ is constant.

We strengthen Theorem~\ref{thm:miracle} in order to also
relate $(\textit{Arith}^{(ab^\flat c)})$ to conditions involving
lattices closed under $f^{-1}$.

A part of the proof of this strengthening has to use either
an implication from Theorem~A or an implication from
Theorem~\ref{th:LatticeGen} which is proved later in the
paper in the framework of general algebras.
For self-containment, we use the implication from
Theorem~\ref{th:LatticeGen}.
Observing that the proof of Theorem~\ref{th:LatticeGen} 
makes no use neither of Theorem~A nor of
Theorem~\ref{thm:stable preserv and lattices} which we
are going to prove, we see that there is no vicious circle
and the proofs of Theorem~\ref{th:LatticeGen} and
Theorem~\ref{thm:stable preserv and lattices} are indeed
self-contained.
\begin{notation}\label{notation Latt emptyset N}
$\Latt_{\suc}^{\emptyset,\N}(L)$ denotes the smallest
sublattice of $\+P(\N)$ which is closed under $\suc^{-1}$
and contains $L\subseteq\N$ and $\emptyset,\N$ 
(i.e. is standardly bounded in the sense of 
Definition~\ref{def standardly bounded lattice}).
\end{notation}
\begin{theorem}\label{thm:stable preserv and lattices}
Let $\+N$ be $\langle \N; \suc \rangle$, 
$\langle \N; + \rangle$ or 
$\langle \N; +, \times \rangle$.
The following conditions are equivalent for any function
$f \colon \N \longrightarrow \N$,
\begin{equation*}
\begin{array}{lll}
(1)&(\textit{PreorderPres})_{\+N}
&
\text{$f$ is $\+N$-stable-preorder preserving}
\\
(2)&(\textit{PreorderPres})_{\+N}^{\finind}
&
\text{$f$ preserves every finite index $\+N$-stable-preorder}
\\
(3)&(\textit{CongPres}^{\uparrow})_{\+N}& 
\left\{\text{\begin{tabular}{l}
$f$ is monotone non-decreasing \\
and $\+N$-congruence preserving
\end{tabular}}\right.
\\
(4)&(\Latt_{\suc}(\forall L))^{\textit{rec},\flat}_{\N}
&\left\{\textit{\begin{tabular}{l}
$f$ is constant or for all recognizable $L \subseteq \N$ \\
$\Latt_{\suc}(L)$ is closed under $f^{-1}$
\end{tabular}}\right.
\\
(5)&{(\Latt^{\emptyset,\N}_{\suc}(\forall L))}^{\textit{rec}}_{\N}
&\textit{$\Latt_{\suc}^{\emptyset,\N}(L)$ is closed under
$f^{-1}$ for all recognizable $L \subseteq \N$}
\\
(6)&(\textit{Arith}^{(ab^\flat c)})
&\left\{\begin{array}{ll}
(a) &\textit{$|y-x|$ divides $|f(y)-f(x)|$ for all
 $x, y \in \N$,}\\
(b^\flat) &\textit{$f$ is constant or $f(x) \geq x$ for all
 $x \in \N$}\\
(c) &\textit{$f$ is monotone non-decreasing}
\end{array}\right.
\end{array}
\end{equation*}
\end{theorem}
\begin{remark}
One can show that there are countably many $\+N$-stable
preorders with infinite index, so that condition $(2)$ is a
non-trivial variant of $(1)$.
On the contrary, the sole $\+N$-congruence with infinite
index is equality on $\N$, so that condition $(3)$ and its
variant dealing with the sole finite index $\+N$-congruences
are trivially equivalent.
\end{remark}
\begin{proof}
Due to Lemma~\ref{l folklore suc + +times idem for rec},
Lemma~\ref{l folklore suc + +times idem}, it suffices to
consider the case $\+N = \langle \N, \suc \rangle$.
Such a reduction is a convenient one since then
$\freez^*(\+N)$ is then the family of iterations of the sole
$\suc$ function. 

Theorem~\ref{thm:miracle} insures that conditions (1), (3),
and (6) are equivalent.

Theorem~\ref{th:LatticeGen} insures equivalence of
conditions in the framework of any algebra $\+A$.
Applying it to the algebra $\+N$, it insures that condition (2)
is equivalent to  the condition
$(\Latt^{\emptyset,A}_{\freez^{*}}(\forall L))_{\+N}$ which
expresses that $f^{-1}(L) \in \Latt_{\+N}^{\emptyset,\N}(L)$
for every $\+N$-recognizable $L \subseteq \N$.
Now, $(\Latt^{\emptyset,A}_{\freez^{*}}(\forall L))_{\+N}$ is
equivalent to condition (5) since we reduced to
$\+N = \langle \N, \suc \rangle$ and closure of a lattice
under $\suc^{-1}$  yields closure under its iterations 
(which constitute $\freez^{*}(\suc)$).
Thus, conditions (2), (5) are equivalent.

Having the equivalences
$(1)\Leftrightarrow(3)\Leftrightarrow(6)$ and
$(2)\Leftrightarrow(5)$, it suffices to prove the equivalences
$(1)\Leftrightarrow(2)$ and $(4)\Leftrightarrow(5)$.

\textbullet\
The implication $(1)\Rightarrow(2)$ is trivial.

\textbullet\
To prove $(2)\Rightarrow(1)$ assume $(2)$ and let
$\sqsubseteq$ be an $\+N$- stable preorder, and suppose
$x \sqsubseteq y$.
Set $n \geq \max(f(x), f(y))$ and define $\sqsubseteq'$ such
that
$\sqsubseteq' \, = \, \sqsubseteq\cup\{(z,t) \mid t \geq n\}$.
It is clearly reflexive and stable and $x \sqsubseteq' y$.
It is also transitive: 
suppose $z \sqsubseteq' t \sqsubseteq' u$,
\\- if $u < n$ then  $z \sqsubseteq t \sqsubseteq u$ 
hence $z \sqsubseteq u$ and therefore $z \sqsubseteq' u$,
\\- if $u \geq n$ then $z \sqsubseteq' u$ by definition.
\\
Finally, $\sqsubseteq'$ has finite index bounded by $n+1$.
So, using $(2)$,
since $x \sqsubseteq' y$, 
we get $f(x) \sqsubseteq' f(y)$ 
hence $f(x) \sqsubseteq f(y)$ because $f(x), f(y) <n$.
This proves $(1)$.

\textbullet\
Let's get $(4)\Rightarrow(5)$. Assume $(4)$.
If $f$ is constant with value $c$ then $f^{-1}(L)$ is either
$\emptyset$ or $\N$ (depending if $c \notin L$ or
$c \in L$) hence
$f^{-1}(L) \in \Latt_{\suc}^{\emptyset,\N}(L)$.
If $f$ is not constant and $L \notin \{\emptyset, \N\}$ then
$(4)$ insures that $f^{-1}$ is in $\Latt_{\suc}(L)$ hence also
in $\Latt_{\suc}^{\emptyset, \N}(L)$.
The case $L$ is $\emptyset$ or $\N$ is trivial.

\textbullet\
Finally, $(5)\Rightarrow(4)$
is proved through a series of seven claims (for readability).

\begin{claimone}
If $L = A$ or $L = A \cup (B+k\N)$ with $A, B$ finite and
$k \geq 2$ then $\emptyset \in \Latt_{\suc}(L)$.
\end{claimone}
\begin{proof}[Proof of Claim 1]
If $t > \max(A)$ then $\suc^{-t}(L) = \emptyset$.
If $t >\max(A\cup B)$ and $L = A \cup (B+k\N)$ then
$\suc^{-t}(L) = C+k\N$ for some
$C \subsetneq \{0, \ldots, k-1\}$.
If $i \in \{0, \ldots, k-1\}\setminus C$ then, setting
$L'=C+k\N$, we have $i-p \notin L'$ for $p = 0, \ldots, i$
and $k+i-p \notin L'$ for $p = i+1, \ldots, k-1$.
Thus, $\bigcap \{\suc^{-q}(L) \mid q = t, \ldots, t+k-1\}
           = \emptyset \in \Latt_{\suc}(L)$.
\end{proof}
\begin{claimtwo}
If $L = A \cup (B+k\N)$ with $B$ non-empty and $k \geq 1$ 
then $\N \in \Latt_{\suc}(L)$.
\end{claimtwo}
\begin{proof}[Proof of Claim 2]
If $i \in B$ then $0 +k\N \subseteq \suc^{-i}(L)$.
Since $\suc^{-n}(k\N) = k-n+k\N$ for $n = 1, \ldots, k-1$,
we see that
$\bigcup_{n = 0, \ldots, k-1} \suc^{-i-n}(L) = k\N$.
\end{proof}
Since 
$\Latt_{\suc}^{\emptyset, \N}(L)
 = \Latt_{\suc}(L) \cup\{\emptyset,\N\}$,
applying Claim 1 and Claim 2 we get
\begin{claimthree}
$\Latt_{\suc}^{\emptyset, \N}(L) = \Latt_{\suc}(L)$ if
$L = A \cup (B+k\N)$ with $A, B$ finite,
$B \neq \emptyset$, and $k \geq 2$.
\end{claimthree}
We now assume the condition
${\Latt^{\emptyset, \N}_{\suc}(\forall L))}^{\textit{rec}}_{\N}$.
\begin{claimfour}
If $f$ satisfies $(5)$ then $f$ is monotonous non-decreasing.
\end{claimfour}
\begin{proof}[Proof of Claim 4]
Fix some $a \in \N$ and let $L = \{ z \mid z \geq f(a)\}$.
Being a final segment of $\N$, the set $L$ is recognizable.
Observe that the $\suc^{-t}(L)$'s, for $t \in \N$, are the final
segments $X_{b} = \{ z \mid z \geq b\}$ where $b \leq f(a)$.
Condition $(5)$ insures that $f^{-1}(L)$ is in the lattice
$\Latt_{\suc}^{\emptyset, \N}(L)$, i.e. is a $(\cup, \cap)$
combination of $\emptyset$, $\N$ and the $X_{b}$'s.
In particular, it is a final segment.
Since obviously $a \in f^{-1}(L)$ this segment is not empty
and is of the form $f^{-1}(L) = \{z \mid z \geq c\}$ for some
$c \leq a$.
Finally, if $x \geq a$ then $x \geq c$ hence $f(x) \in L$,
which means that $f(x) \geq f(a)$.
\end{proof}
\begin{claimfive}
If $f$ satisfies $(5)$ and is not constant then it takes
infinitely many values.
\end{claimfive}
\begin{proof}[Proof of Claim 5]
By way of contradiction, suppose that the range of $f$ is a
finite set $\{b_{1}, \ldots, b_{n}\}$ with $n \geq 2$ and
$b_{1} < \ldots < b_{n}$.
By monotony of $f$ the inverse image of $\{ b_{n}\}$ is a
final segment of $\N$ and since $n \geq 2$ this is a strict
final segment of $\N$, say $f^{-1}(b_{n}) = p+\N$ with
$p \geq 1$.
Let $L = b_{n} + k\N$ with $k > b_{n}$.
Then $f^{-1}(L) = f^{-1}(b_{n}) = p+\N$.
Observe that the $\suc^{-t}(L)$'s, with $t \in \N$, are the
sets
\begin{align*}
\suc^{-t}(L) =&
\left\{\begin{array}{ll}
b_{n}-t+k\N
&\textit{if $t = 0, \ldots, b_{n}-1$,}
\\
((b_{n}-t) \bmod k)+k\N
&\textit{if $t \geq b_{n}$.}
\end{array}\right.
\end{align*}
Thus, the $\suc^{-t}(L)$'s are all the sets $a+k\N$ with
$a \leq \max(b_{n}, k-1) = k-1$ since we choose $k > b_{n}$.
Considering intersections of $\emptyset$, $\N$, and the
$\suc^{-t}(L)$'s does not produce any set outside these
ones.
Since $(\cup, \cap)$-combinations can be normalized to
finite non-empty unions of finite non-empty intersections,
we see that sets in the lattice 
$\Latt_{\suc}^{\emptyset, \N}(L)$ are finite non-empty
unions of $\emptyset$, $\N$, and the $\suc^{-t}(L)$'s. 
If $f^{-1}(L) = p+\N$ were in the lattice 
$\Latt_{\suc}^{\emptyset, \N}(L)$, it would be a finite
non-empty union of sets $a+k\N$ with $a \leq k-1$ since
$\emptyset$ is useless in such a union and $\N$ would
give a strict superset of $p+\N$ because $p \geq 1$.
Thus, we would have $f^{-1}(L) = \bigcup_{a \in I}a+k\N$
with $I \subseteq \{0, \ldots, k-1\}$.

Since 
$f^{-1}(L) = p+\N = \bigcup_{i=0, \ldots, k-1 }p+i+k\N$
with $p \geq 1$, we would have
\begin{equation}\label{eq if f not ct}
\bigcup_{a \in I} a+k\N
= \bigcup_{i= 0, \ldots, k-1}p+i+k\N
\quad\textit{where $I \subseteq \{0, \ldots, k-1\}.$}
\end{equation}
The right side is a union of $k$ many pairwise disjoint sets
$p+i+k\N$ for $i = 0, \ldots, k-1$.
Equality~\eqref{eq if f not ct} implies that
$I = \{p, \ldots, p+k-1\}$.
Contradiction since $I \subseteq \{0, \ldots, k-1\}$ and
$p \geq 1$.
\end{proof}
\begin{claimsix}
If $f$ satisfies $(5)$ and $L$ is finite then $f^{-1}(L)$ is
finite and $f^{-1}(L) \in \Latt_{\suc}(L)$.
\end{claimsix}
\begin{proof}[Proof of Claim 6]
Since $f$ takes infinitely many values and is monotone, the
inverse images of singletons are finite (possibly empty)
segments of $\N$.
In particular, $f^{-1}(L)$ is finite.
Assumption $(5)$ insures that $f^{-1}(L)$ is a finite
non-empty union of finite non-empty intersections of
$\emptyset$, $\N$, and the $\suc^{-t}(L)$'s with $s \in \N$.
Since $f^{-1}(L)$ is finite we can exclude $\N$ from this
combination, and using Claim 1 we can also exclude
$\emptyset$.
Thus, $f^{-1}(L)$ is in fact in $\Latt_{\suc}(L)$.
\end{proof}
\begin{claimseven}
If $f$ satisfies $(5)$ and if $L = A \cup (B+\N)$ with $A, B$
finite and $B$ not empty then $f^{-1}(L)$ is cofinite and
$f^{-1}(L) \in \Latt_{\suc}(L)$.
\end{claimseven}
\begin{proof}[Proof of Claim 7]
Using assumption $(5)$, consider a finite non-empty union
of finite non-empty intersections of $\emptyset$, $\N$, and
the $\suc^{-t}(L)$'s with $t \in \N$ which gives $f^{-1}(L)$.
Claim 2 insures that $\N \in \Latt_{\suc}(L)$ hence we can
exclude $\N$ and only use $\emptyset$ and the
$\suc^{-t}(L)$'s.
Now, since $L$ is cofinite, $\N \setminus L$ is finite and by
Claim 6 $f^{-1}(\N \setminus L)$ is also finite.
Hence $f^{-1}(L)$ is cofinite.
In particular, it is not empty and we can exclude 
$\emptyset$ from the combination giving $f^{-1}(L)$.
This shows that $f^{-1}(L) \in \Latt_{\suc}(L)$.
\end{proof}
By Lemma~\ref{l rec N} and Claims 3, 6, 7, in all cases
$f^{-1}(L) \in \Latt_{\suc}(L)$, proving $(4)$.
\end{proof}

\section{Lattice and Boolean algebra conditions for stable
preorder and congruence preservation in an algebra} 
\label{s lattice cong preserv}
%
\subsection{Lattices $\Latt_{\+A}^{\emptyset,A}(L)$,
                              $\Latt_{\+A}^{\infty}(L)$, 
Boolean algebras $\Bool_{\+A}^{\emptyset,A}(L)$, 
                             $\Bool_{\+A}^{\infty}(L)$}
\label{ss lattices of subsets close by preimage}
%
As witnessed by the lattice of closed subsets of a topological
space such as $\RR$, the supremum of an infinite family in a
lattice of sets may strictly contain its set theoretic union. This
is why we recall the following notion.
\begin{definition}
\label{def standardly bounded set complete}
1. A {\em lattice} $\Latt$ of subsets of $A$ is
{\em set-complete} if for every (possibly empty or infinite)
$\+F\subseteq\Latt$ the union $\bigcup_{X \in \+F}X$ and
the intersection $\bigcap_{X \in \+F}X$ are in $\Latt$
(hence are the supremum and infimum of $\+F$).
In particular, such a lattice is standardly bounded
(cf. Definition~\ref{def standardly bounded lattice}).

2. A standardly bounded lattice closed under
complementation is called a {\em standardly bounded
Boolean algebra}.

3. A set-complete (hence standardly bounded) lattice 
closed under complementation is called a {\em set
complete Boolean algebra}.
\end{definition}
Quitting the rich framework of algebras
$\langle \N, \suc \rangle$, $\langle \N, + \rangle$,
$\langle \N, +, \times \rangle$, and considering general
algebras, a key tool is the notion of freezification
introduced in \S\ref{ss reducing to arity 1}.
\begin{definition}\label{def L B}
For $\+A = \langle A; \Xi \rangle$ an algebra and
$L \subseteq A$, among the families of subsets of $A$
which contain $L$ and are closed under all the
$\gamma^{-1}$'s for $\gamma \in \freez^{*}(\+A)$, we
respectively denote 
\\\centerline{
\begin{tabular}{ll}
$\Latt_{\+A}(L)$ & the smallest lattice,
\\
$\Latt_{\+A}^{\emptyset, A}(L)$
& the smallest standardly bounded lattice
(i.e. $\Latt_{\+A}(L) \cup \{\emptyset, A\}$),
\\
 $\Latt_{\+A}^{\infty}(L)$ 
& the smallest set complete (hence standardly bounded)
lattice,
\\
$\Bool_{\+A}^{\emptyset, A}(L)$
& the smallest standardly bounded Boolean algebra,
\\
$\Bool_{\+A}^{\infty}(L)$
& the smallest set complete (hence standardly bounded)
Boolean algebra.
\end{tabular}}
\end{definition}
The following Lemma is straightforward.
\begin{lemma}\label{l if finite all equal}
If $\Latt_{\+A}(L)$ is finite then
$$\Latt_{\+A}^{\emptyset, A}(L)
= \Latt_{\+A}^{\infty}(L)
\text{\quad and\quad }
\Bool_{\+A}^{\emptyset, A}(L)
= \Bool_{\+A}^{\infty}(L).$$
\end{lemma}
\begin{figure}[h]
\center
\fbox{$\begin{array}{crrcll}
\langle \N; \suc \rangle, \langle \N; +\rangle
&X\mapsto &X-n&=&\{y \mid n+y \in X\}&\textit{with }n
 \in \N
\\
\langle\N; \times \rangle&X \mapsto &X/n&=&\{ y \mid ny
 \in X\}
&\textit{with }n \in \N
\\
\langle \N; +, \times \rangle&X \mapsto &(X-m)/n&=&
\{y \mid ny+m \in X\}&\textit{with }m, n \in \N
\\
\langle \Sigma^*; \cdot \rangle&X \mapsto &u^{-1}
 X v^{-1}&=
& \{x \mid uxv \in X\}&\textit{with }u, v \in \Sigma^*
\end{array}$}
\caption{Maps $X \mapsto \gamma^{-1}(X)$ 
for $\gamma \in \freez^*(\+A)$
\label{gamma-1a}}
\end{figure}
\begin{note}
Maps $X \in \+P(A) \mapsto \gamma^{-1}(X)$ are called
{\em cancellations} in \cite{almeida}, \S3.1, p.55.
See Figure~\ref{gamma-1a} for examples.
\end{note}
\begin{lemma}\label{bool:saturated}
Let $\sim$ be an $\+A$-congruence.
If $L\subseteq\+A$ is $\sim$-saturated then so is every set
in $\Bool_{\+A}^\infty(L)$.
In particular, this is true for the syntactic congruence
$\sim_{L}$.
\end{lemma}
\begin{proof}
The family of $\sim$-saturated subsets of $A$ is closed
under arbitrary unions, complementation and preimage by
maps in $\freez^{*}(\+A)$, hence it contains
$\Bool_{\+A}^\infty(L)$.
\end{proof}
\begin{lemma}\label{preorder:initial}
Let $\+A$ be an algebra and $\preceq$ an $\+A$-stable
preorder. 
If $L$ is a $\preceq$-initial segment (i.e. $a \in L$ and
$b \preceq a$ imply $b \in L$) then so is every set in
$\Latt_{\+A}^\infty(L)$.
\end{lemma}
\begin{proof}
The family of $\preceq$-initial segments of $A$ is closed
under arbitrary unions and intersections and preimage by
maps in $\freez^{*}(\+A)$, hence it contains
$\Latt_{\+A}^\infty(L)$.
\end{proof}
\begin{lemma}\label{l;recFinite}
If $L$ is $\+A$-recognizable then $\Bool_\+A(L)$ is finite, 
(so that equalities of Lemma~\ref{l if finite all equal} hold)
and every set in $\Bool_\+A(L)$ is also $\+A$-recognizable.
\end{lemma}
\begin{proof}
By Proposition~\ref{p:syntactic largest}, $L$ is
$\sim_L$-saturated and since $L$ is $\+A$-recognizable the
congruence $\sim_L$ has finite index, say $k$, and the
$2^k$ many $\sim_L$-saturated subsets of $A$ are all
$\+A$-recognizable.
By Lemma~\ref{bool:saturated}, all sets in
$\Bool^\infty_\+A(L)$ are $\sim_L$-saturated, hence this
algebra is finite with cardinal at most $2^k$, constituted
of $\+A$-recognizable sets.
\end{proof}

%
\subsection{Representing $f^{-1}(L)$ when $f$ preserves
$\sim_{L}$ or $\preceq_L$}
\label{ss f-1L preceq}
%
\begin{lemma}\label{l f^-1(L)}
Consider an algebra $\+A = \langle A; \Xi \rangle$, a
function $f \colon A \to A$ and a set $L \subseteq A$.  
Let $X_{a,L} = \{\gamma \in \freez^{*}(\Xi) \mid a \in
 \gamma^{-1}(L)\}$ and
$Y_{a,L} = \{\gamma \in \freez^{*}(\Xi) \mid a \notin
  \gamma^{-1}(L)\}$.

1. If $f$ preserves the syntactic congruence $\sim_L$ then 
\begin{align}\label{eq f-1L sim}
f^{-1}(L) &= \bigcup_{a \in f^{-1}(L)} \{x \mid x \sim_L a\} 
	 	= \bigcup_{a \in f^{-1}(L)} \
  \left(   \bigcap_{\gamma \in X_{a,L}} \gamma^{-1}(L)
  \cap
  \bigcap_{\gamma \in Y_{a,L}}
    (A \setminus \gamma^{-1}(L))
  \right)
\end{align}
In particular, $f^{-1}(L) \in \Bool_{\+A}^{\infty}(L)$. 

Moreover, if $L$ is $\+A$-recognizable then 
$\Bool_{\+A}^{\emptyset,A}(L) =
 \Bool_{\+A}^{\infty}(L)$
and is finite.

2. If $f$ preserves the syntactic preorder $\preceq_L$ then
$f^{-1}(L) \in \Latt_{\+A}^\infty(L)$ and
\begin{align}\label{eq f-1L preceq}
f^{-1}(L) &= \bigcup_{a \in f^{-1}(L)} \{x \mid x \preceq_L a\}
                = \bigcup_{a \in f^{-1}(L)} \ 
                   \bigcap_{\gamma \in X_{a,L}} \gamma^{-1}(L)
\end{align}
In particular, $f^{-1}(L) \in \Latt_{\+A}^{\infty}(L)$. 

Moreover, if $L$ is $\+A$-recognizable then 
$\Latt_{\+A}^{\emptyset, A}(L) =
 \Latt_{\+A}^{\infty}(L)$
and is finite.
\end{lemma}
\begin{proof} 
1. If $f^{-1}(L) = \emptyset$ then
$\bigcup_{a \in f^{-1}(L)} \{x \mid x \sim_L a\}$ is empty as
union of an empty family, hence is equal to $f^{-1}(L)$.
Assume now $f^{-1}(L) \neq \emptyset$.
Since $a \sim_L a$ we have 
$f^{-1}(L) \subseteq \bigcup_{a \in f^{-1}(L)} \{x \mid
 x \sim_{L} a\}$.
Conversely, if $x \sim_{L} a$ where $a \in f^{-1}(L)$ then, as
$f$ preserves $\sim_L$, we also have $f(x) \sim_L f(a)$, which
means that, for every $\gamma \in \freez^*(\+A)$, we have
the equivalence
$\gamma(f(a)) \in L \Leftrightarrow \gamma(f(x)) \in L$. 
For $\gamma = id$, since $f(a) \in L$ this gives $f(x) \in L$,
i.e. $x \in f^{-1}(L)$. 
Thus, $\{x \mid x \sim_{L} a\} \subseteq f^{-1}(L)$ for every
$a \in f^{-1}(L)$, which proves 
$\bigcup_{a \in f^{-1}(L)} \{x \mid x \sim_{L} a\}
 \subseteq f^{-1}(L)$,
giving the first equality in \eqref{eq f-1L sim}.
The second equality simply translates the definition of
$\sim_{L}$ (cf. Definition~\ref{def syntactic cong pre}).
 
Finally, the second equality in \eqref{eq f-1L sim} yields
$f^{-1}(L) \in \Bool_{\+A}^{\infty}(L)$, the case of empty
union (in case $f^{-1}(L)$ is empty) and empty
intersections (in case some $X_{a,L}$'s or $Y_{a,L}$'s are
empty) justifying the requirement that $\emptyset$ and
$A$ be forcibly put in this Boolean algebra.

For item 2 argue similarly and for the case of recognizable
$L$ use Lemma~\ref{l;recFinite}.
\end{proof}

\subsection{$\+A$-congruence (resp. $\+A$-stable
preorder) preservation and Boo\-le\-an algebras (resp.
lattices) of subsets of $\+A$}
\label{ss:preserv and latt BA}
%
\begin{table}[h]
\begin{center}
$\begin{array}{l|ll|}
&&\text{\quad Four pairs of equivalent conditions}
\\
\cline{2-3}
\text{\scriptsize(1)}&
(\textit{PreordPres})_{\+A}
&\textit{$f$ is $\+A$-stable-preorder preserving}
\\
&(\Latt^{\infty}_{\freez^{*}}(\forall L))_{\+A}
& \textit{$f^{-1}(L) \in
 \Latt_{\+A}^{\infty}(L)$
          for every $L \subseteq A$}
\smallskip
\\
\cline{2-3}
\text{\scriptsize(2)}&
(\textit{PreordPres})^{\finind}_{\+A}
&\textit{$f$ preserves all finite index
 $\+A$-stable-preorders}
\\
&(\Latt^{\emptyset,A}_{\freez^{*}}(\forall L))_{\+A}^{\textit{rec}}
& \textit{$f^{-1}(L) \in \Latt_{\+A}^{\emptyset,A}(L)$
          for every recognizable $L \subseteq A$}
\smallskip
\\
\cline{2-3}
\text{\scriptsize(3)}&
(\textit{CongPres})_{\+A}
&\textit{$f$ is $\+A$-congruence preserving}
\\
&(\Bool^{\infty}_{\freez^{*}}(\forall L))_{\+A}
& \textit{$f^{-1}(L) \in
 \Bool_{\+A}^{\infty}(L)$
          for every $L \subseteq A$}
\smallskip
\\
\cline{2-3}
\text{\scriptsize(4)}&
(\textit{CongPres})^{\finind}_{\+A}
&\textit{$f$ preserves all finite index $\+A$-congruences}
\\
& (\Bool^{\emptyset,A}_{\freez^{*}}(\forall L))_{\+A}^{\textit{rec}}
& \textit{$f^{-1}(L) \in \Bool_{\+A}^{\emptyset, A}(L)$
          for every recognizable $L \subseteq A$}
\\\cline{2-3}
\end{array}$

\medskip

\noindent
\begin{tikzpicture}\label{figure equivalences}
\coordinate (P1) at (0,6) ;		
\coordinate (P1C1) at (0,5.7) ;	
\coordinate (P1Q1) at (0.3,5.7) ;
\coordinate (P1P2) at (1.3,6) ;	

\coordinate (P2) at (7.4,6) ;		
\coordinate (P2C2) at (7.4,5.7) ;	
\coordinate (P2Q2) at (7.1,5.7) ;
\coordinate (P2P1) at (5.8,6) ;	

\coordinate (C1) at (0,0) ;		
\coordinate (C1P1) at (0,0.3) ; 	
\coordinate (C1D1) at (0.3,0.3) ; 
\coordinate (C1C2) at (1.2,0) ;	

\coordinate (C2) at (7.4,0) ;	
\coordinate (C2C1) at (5.8,0) ;	
\coordinate (C2D2) at (7.1,0.3) ; 
\coordinate (C2P2) at (7.4,0.3) ; 

\coordinate (Q1) at (1.5,4) ; 
\coordinate (Q1Q2) at (2.7,4) ; 
\coordinate (Q1D1) at (1.6,3.7) ;
\coordinate (Q1P1) at (1.4,4.3) ; 

\coordinate (Q2) at (5.8,4) ; 
\coordinate (Q2Q1) at (4.5,4) ; 
\coordinate (Q2D2) at (5.7,3.7) ; 
\coordinate (Q2P2) at (5.8,4.3) ;

\coordinate (D1) at (1.5,2) ; 
\coordinate (D1D2) at (2.7,2) ;
\coordinate (D1Q1) at (1.6,2.3) ;
\coordinate (D1C1) at (1.7,1.7) ;

\coordinate (D2) at (5.8,2) ; 
\coordinate (D2D1) at (4.5,2) ; 
\coordinate (D2Q2) at (5.7,2.3) ; 
\coordinate (D2C2) at (5.8,1.7) ; 

\draw (P1) node{$(\textit{PreordPres})_{\+A}$} ;
\draw (P2) node{$(\textit{PreordPres})^{\finind}_{\+A}$} ;
\draw (C1) node{$(\textit{CongPres})_{\+A}$} ;
\draw (C2) node{$(\textit{CongPres})^{\finind}_{\+A}$} ;
\draw (Q1) node{$(\Latt^{\infty}_{\freez^{*}}(\forall L))_{\+A}$} ;
\draw (Q2) node{$(\Latt^{\emptyset, A}_{\freez^{*}}(\forall L))_{\+A}^{\textit{rec}}$} ;
\draw (D1) node{$(\Bool^{\infty}_{\freez^{*}}(\forall L))_{\+A}$} ;
\draw (D2) node{$(\Bool^{\emptyset, A}_{\freez^{*}}(\forall L))_{\+A}^{\textit{rec}}$} ;

\draw [thick,double,->] (P1P2) -- (P2P1) ; 
\draw [thick,double,->] (Q1Q2) -- (Q2Q1) ;
\draw [thick,double,->] (D1D2) -- (D2D1) ;
\draw [thick,double,->] (C1C2) -- (C2C1) ;
\draw [thick,double,->] (P1C1) -- (C1P1) ; 
\draw [thick,double,->] (Q1D1) -- (D1Q1) ;
\draw [thick,double,->] (Q2D2) -- (D2Q2) ;
\draw [thick,double,->] (P2C2) -- (C2P2) ;
\draw [thick,double,<->] (P1Q1) -- (Q1P1) ; 
\draw [thick,double,<->] (P2Q2) -- (Q2P2) ;
\draw [thick,double,<->] (C1D1) -- (D1C1) ;
\draw [thick,double,<->] (C2D2) -- (D2C2) ;
\end{tikzpicture}
\caption{Implications and equivalent conditions of
Theorem~\ref{th:LatticeGen}}
\label{table equivalences}
\end{center}
\end{table}
\begin{theorem}[Preservation in general algebras]
\label{th:LatticeGen} 
Let $\+A = \langle A; \Xi \rangle$ be an algebra and
$f \colon A \to A$.
For each one of the pairs of conditions detailed in
Table~\ref{table equivalences}, the two shown conditions
are equivalent.
The figure of Table~\ref{table equivalences} illustrates these
equivalences and the straightforward implications.
\end{theorem}
\begin{proof} 
All the simple implications (horizontal and vertical arrows)
of the diagram are trivial.
So, we are left with the four inclined bi-implications to prove.
\\
\textbullet\ (1)\quad
$(\textit{CongPres})_{\+A}\Longrightarrow
 (\Bool^{\infty}_{\freez^{*}}(\forall L))_{\+A}$.
Assume $f$ preserves all $\+A$-congruences. 
Then $f$ preserves the syntactic congruence $\sim_{L}$ of
any $L$ and item 1 in Lemma \ref{l f^-1(L)}  insures that
$f^{-1}(L)\in\Bool^{\infty}_{\+A}$.

$(\Bool^{\infty}_{\freez^{*}}(\forall L))_{\+A}
 \Longrightarrow	(\textit{CongPres})_{\+A}$.
Let $L = \{x \mid x \sim f(a)\}$ where $\sim$ is an
$\+A$-congruence and $a \in A$.
Condition $(\Bool^{\infty}_{\freez^{*}}(\forall L))_{\+A}$
insures that $f^{-1}(L)$ is in $\Bool_{\+A}^{\infty}(L)$.
Since $L$ is $\sim$-saturated so is also every set in
$\Bool_{\+A}^{\infty}(L)$ by Lemma~\ref{bool:saturated}.
In particular $f^{-1}(L)$ is $\sim$-saturated.
Now, $a \in f^{-1}(L)$ since $f(a) \sim f(a)$, hence if
$x \sim a$ then $x$ is also in $f^{-1}(L)$, hence $f(x) \in L$
which means $f(x) \sim f(a)$.
Thus, $f$ preserves the congruence $\sim$.
\\
\textbullet\ (2)\quad
$(\textit{CongPres})^{\finind}_{\+A}\Longleftrightarrow
	(\Bool^{\emptyset,A}_{\freez^{*}}(\forall L))_{\+A}^{\textit{rec}}$.
Since $L$ is recognizable if and only if $\sim_{L}$ has finite
index (cf. Proposition~\ref{p:syntactic largest}), we can argue
as in (1) above.
\\
\textbullet\ (3)\
$(\textit{PreordPres})_{\+A} \Longrightarrow
 (\Latt_{\freez^{*}}^{\infty}(\forall L))_{\+A}$.
Argue as for congruences and Boolean algebras, using
item 2 in Lemma \ref{l f^-1(L)}.

$(\Latt_{\freez^{*}}^{\infty}(\forall L))_{\+A}
 \Longrightarrow	(\textit{PreordPres})_{\+A}$.
Let $\preceq$ be an $\+A$-stable preorder and $a \in A$
and $L = \{x \mid x \preceq f(a)\}$.
Condition $(\Latt_{\freez^{*}}^{\infty}(\forall L))_{\+A}$
insures that $f^{-1}(L)$ is in $\Latt_{\+A}^{ \infty}(L)$.
Since $L$ is a $\preceq$-initial segment so is also every set
in $\Latt_{\+A}^{\emptyset, A, \infty}(L)$ (by
Lemma~\ref{preorder:initial}).
In particular $f^{-1}(L)$ is a $\preceq$-initial segment.
Now, $a \in f^{-1}(L)$ since $f(a) \preceq f(a)$, hence if
$x \preceq a$ then $x$ is also in $f^{-1}(L)$, hence
$f(x) \in L$ which means $f(x) \preceq f(a)$.
Thus, $f$ preserves the stable preorder $\preceq$.
\\
\textbullet\ (4)\quad
$(\textit{PreordPres})^{\finind}_{\+A} \Longleftrightarrow
	(\Latt_{\freez^{*}}^{\emptyset, A}(\forall L))_{\+A}^{\textit{rec}}$.
Since $L$ is recognizable if and only if $\preceq_{L}$ has
finite index (cf. Proposition~\ref{p:syntactic largest}), we can
argue as in (3) above.
\end{proof}
\begin{remark}\label{rk if f constant}
In case $f$ is constant it obviously preserves congruences
and stable preorders.
If its value is $a$ then the inverse image $f^{-1}(L)$ is
either the empty set if $a \notin L$ or the whole set $A$ if
$a \in L$.
This illustrates why we have to put $\emptyset$ and $A$ in
the lattices and Boolean algebras considered in
Theorem~\ref{th:LatticeGen}.
\end{remark}
%

%
\section{Adding algebraic hypothesis to enrich
Table~\ref{table equivalences}}
\label{s improving}
%
In this section we consider particular properties of algebras 
allowing to replace some implications in the figure of 
Table~\ref{table equivalences} by logical equivalences.
%
\subsection{Algebras with a group operation and the right
side-trapezium}
\label{ss:when there is a group operation}
%
For some algebras it happens that the notion of finite index
stable preorder is equivalent to that of finite index
congruence, cf.
Propositions~\ref{p:stable cancellable preorder with finite index}
and \ref{p unit ring} below.

We first give a version involving semigroups and cancellable
stable preorders.
\begin{definition}\label{def:cancellative}
1) A {\em semigroup} $S$ is said to be {\em cancellable} if
$xz = yz$ implies $x = y$ and $zx = zy$ implies $x = y$.

2) A {\em stable preorder} $\preceq$ on $S$ is said to be
{\em cancellable} if $xz \preceq yz$ implies $x \preceq y$
and $zx \preceq zy$ implies $x \preceq y$.
\end{definition}
Let's recall two folklore Lemmas.
\begin{lemma}\label{l:stable order on finite group}
The only stable order of a finite group $G$ is the identity
relation.
\end{lemma} 
\begin{proof}
Assuming $x \preceq y$ we prove $x = y.$ 
Let $e$ be the unit of $G$. 
Stability under left product by $x^{-1}$ and $(x^{-1} y)^{n}$
successively yield $e \preceq x^{-1}y$ and then
$(x^{-1} y)^{n} \preceq (x^{-1} y)^{n+1}$.
By transitivity, $e \preceq x^{-1} y \preceq (x^{-1} y)^{n}$.
As the group is finite there exists $k$ and $n \geq 1$ such
that $(x^{-1} y)^{k+n} = (x^{-1} y)^{k}$ hence
$(x^{-1} y)^{n} = e$.
Thus, $e \preceq x^{-1}y \preceq e$ and by antisymmetry
$e = x^{-1} y$ and $x = y$.
\end{proof} 
\begin{lemma}\label{l:finite cancellative semigroup}
Any finite cancellable semigroup $S$ is a group.
\end{lemma} 
\if 42
\begin{proof} 
Cancellability ensures that, for every $a\in S$,
the maps $x\mapsto ax$ and $x\mapsto xa$ are injective 
hence are bijections $S \to S$ because $S$ is finite. 
In particular, for all $a\in S$ 
there exist $e'_{a}$ and $e''_{a}$ such that $ae'_{a}=a$ and $e''_{a}a=a$.
For $a,b\in S$ we then have $ae'_{a}b = ab = ae''_{b}b$
and by cancellability $e'_{a} = e''_{b}$.
Letting $b$ be $a$ we see that $e'_{a} = e''_{a}$, 
and fixing $a$ and varying $b$, we see that $e'_{a}$ does not depend on $a$. 
This proves that the common value of the $e'_{a}$ is a unit $e$ of $S$.
Also, for $a\in S$ there exists $a'$ and $a''$ such that $aa'=a''a=e$.
Then $a'=ea'=(a''a)a'=a''(aa')=a''e=a''$, proving that $a'$ is an inverse of $a$.
\end{proof}
\fi

\begin{proposition}
\label{p:stable cancellable preorder with finite index}
Let $\+A = \langle A; \Xi\rangle$ be an algebra where $\Xi$
contains a semigroup operation.
Every cancellable finite index stable preorder $\preceq$ of
$\+A$ coincides with its associated congruence $\sim$.
\end{proposition}
\begin{proof}
The semigroup operation on $\+A$ induces a semigroup
operation on the quotient algebra $\+G = \+A/{\sim}$ with
carrier set $G = A/{\sim}$.

The cancellability property of the congruence $\sim$ 
implies that of the semigroup operation on $G$.
Indeed, suppose $X, Y, Z \in G$ satisfy $XZ = YZ$ and let
$x, y, z \in A$ be representatives of the classes $X, Y, Z$. 
Then we have $xz \sim yz$ and cancellability yields
$x \sim y$, hence $X = Y$. Idem if $ZX = ZY$.

As $\sim$ has finite index, $G$ is finite and
Lemma~\ref{l:finite cancellative semigroup} insures that the
semigroup operation on $G$ is a group operation.
As $\sim$ is the congruence associated to the stable
preorder $\preceq$, it induces a quotient stable order
$\preceq/\!\!\sim$ on the finite quotient algebra $\+G$.
As $\+G$ is an expansion of a finite group,
Lemma~\ref{l:stable order on finite group} ensures that
$\preceq/\!\!\sim$ is the identity relation on $G$, hence
$\preceq$ coincides with $\sim$ and is therefore a
congruence on $\+A$.
\end{proof}
In the framework of groups the requirement that a preorder
be stable is quite a strong requirement as shown by the
following corollary of 
Proposition~\ref{p:stable cancellable preorder with finite index}.
\begin{proposition}
\label{p:stable finite index preorder in group}
Let $\+A = \langle A; \Xi \rangle$ be an algebra where $\Xi$
contains a group operation.

1. Every stable preorder with finite index is a congruence.

2. A function $f \colon A\to A$ is $\+A$-stable finite index
preorder preserving if and only if it is finite index
$\+A$-congruence preserving.
In particular, the right side trapezium in
Table~\ref{table equivalences} becomes a series of logical
equivalences.

3. If $L \subseteq A$ is recognizable then its syntactic
preorder is equal to its syntactic congruence and 
$\Latt_{\+A}^{\emptyset, A}(L) = \Latt_{\+A}^{\infty}(L)$ and
$\Bool_{\+A}^{\emptyset, A}(L) = \Bool_{\+A}^{\infty}(L)$.
\end{proposition}
\begin{proof}
1. Observe that stability of $\preceq$ implies its
cancellability: if $xz \preceq yz$ then
$xzz^{-1} \preceq yzz^{-1}$ hence  $x \preceq y$.
Idem if $zx \preceq zy$. Then apply
Proposition~\ref{p:stable cancellable preorder with finite index}.

2. Obvious corollary of 1.

Finally, for 3 observe that the syntactic congruence of a
recognizable set has finite index.
\end{proof}
%

%
\subsection{Algebras expanding unit rings and the two
side-trapeziums}
\label{ss:when there is a unit ring operation}
%
The following easy folklore result is a kind of companion to 
Proposition~\ref{p:stable finite index preorder in group}.
\begin{proposition}\label{p unit ring}
If the algebra $\+A$ is an algebraic expansion of a unit ring
then every $\+A$-stable preorder is an $\+A$-congruence.
As a consequence, in the two side trapeziums in
Table~\ref{table equivalences} all arrows become logical
equivalences.
\end{proposition}
\begin{proof}
Let $1$ be the unit of the unit ring multiplication and $-1$
its opposite relative to the ring addition.
If $\preceq$ is an $\+A$-stable preorder and $x \preceq y$
then $(-1)x \preceq (-1)y$ hence
$y = x+(-1)x+y \preceq x+(-1)y+y = x$.
Thus, $x \preceq y$ and $y \preceq x$, showing that
$\preceq$ is an $\+A$-congruence.
\end{proof}

\subsection{Residually finite, residually c-finite and
sp-finite algebras}
\label{ss residually finite}
%
In this section we review some notions of residually finite
algebra stronger than the classical one, as in
\cite{almeida} (p. 23) and \cite{KaarliPixley} (p. 102), that
are used in the next
\S\ref{ss improving for c-residually finite}.
\begin{definition}[Classical definition]
\label{def:residually finite usual}
An algebra $\+A$ is {\em residually finite} if morphisms from
$\+A$ to finite algebras separate points, i.e., if $x \neq y$
then there exists a morphism $\varphi$ into a finite algebra 
such that $\varphi(x) \neq \varphi(y)$.
In other words, $\+A$ is residually finite if equality is the
intersection of congruences having finite index.
\end{definition}
\begin{lemma}\label{ex res fin algebras}
Free groups and free monoids (in particular
$\langle \N; + \rangle$ and $\langle \Z; + \rangle$) are
residually finite whereas the algebra
$\langle \Q; + \rangle$ is not residually finite.
\end{lemma}
\begin{proof}
{\it Case of free groups.}
See Daniel E. Cohen's book \cite{Cohen89}, pp. 7 and 11.

{\it Case of free monoids.}
Let's recall the easy folklore proof.
Given words $u \neq v$ in the set $\Sigma^{*}$ of finite
words over alphabet $\Sigma$, let $\Delta$ be the finite set
of letters of $\Sigma$ used in the words $u, v$ and
$k = \max(|u|, |v|)$ be the maximum of the lengths of $u$
and $v$.
Consider the set $I$ of all words containing a letter in
$\Sigma\setminus\Delta$ or having length strictly greater
than $k$.
This set $I$ is a two-sided ideal of the monoid $\Sigma^{*}$
and the quotient monoid $M = \Sigma^{*}/I$ is finite since
it can be identified to
$\{x \in \Delta^{*} \mid |x| \leq k\} \cup\{0\}$ where $0$ is
a zero element and the product of $x$ and $y$ is $xy$ if
$|xy| \leq k$ and $0$ otherwise.
Clearly, the canonical morphism
$\varphi \colon \Sigma^{*} \to M$ separates $u$ and $v$ 
since it maps the distinct elements $u$ and $v$ onto
themselves. 

{\it Case of $\langle \Q; + \rangle$.}
Any congruence of a group is defined by the class of its
neutral element which is a normal subgroup.
It is known that in $\langle \Q; + \rangle$ the sole finite
index subgroup is $\Q$ itself. 
Thus, the sole finite index congruence is the trivial one
which identifies all elements of $\Q$ hence
$\langle \Q; + \rangle$ is not residually finite.
\end{proof}
Extending Definition~\ref{def:residually finite usual} to all
congruences --~not only equality~-- we get the following
notions.
\begin{definition}\label{def:residually finite c sp}
1) A {\em congruence} $\sim$ on an algebra $\+A$ is said
to be {\em c-residually finite} if it is the intersection of a
family of finite index congruences. 
In other words, a congruence $\sim$ is c-residually finite if
the quotient $\+A/\sim$ is residually finite in the usual
sense of Definition~\ref{def:residually finite usual}.

2) A {\em stable preorder} on an algebra $\+A$ is said to be
{\em sp-residually finite} if it is the intersection of a family of
stable preorders all of which have finite index.

3) An {\em algebra} $\+A$ is said to be {\em c-residually
finite} (resp. {\em sp-residually finite}) if all congruences
(resp. stable preorders)  on $\+A$ are c-residually finite
(resp. sp-residually finite).
\end{definition}
Every congruence being a preorder,
Definition~\ref{def:residually finite c sp} gives a priori two
notions of residual finiteness for a congruence. 
However, both notions coincide.
\begin{lemma}\label{l:p residually finite implies c}
1. If a stable preorder $\preceq$ is sp-residually finite then
its associated congruence $\sim$ is c-residually finite.

2. A congruence is c-residually finite if and only if, as a
preorder, it is sp-residually finite. In particular, an
sp-residually finite algebra is also c-residually finite.
\end{lemma}
\begin{proof}
1) Let $(\preceq_i)_{i \; \in \; I}$ be a family of stable
preorders having finite index and such that
${\preceq} = {\bigcap_{i \; \in \; I} \preceq_i}$.
Let $\sim_i$ be the congruence associated to $\preceq_i$.
We show that ${\sim}={\bigcap_{i \; \in \; I} \sim_i}$.
As $\sim_i$ is included in $\preceq_i$ we have
$(\bigcap_{i \; \in \; I} \sim_i) \subseteq
 (\bigcap_{i \; \in \; I} \preceq_i) = {\preceq}$.
As $\bigcap_{i\in I} \sim_i$ is a congruence, the last
nclusion yields
$(\bigcap_{i \; \in \; I} \sim_i) \subseteq {\sim}$.
The inclusion of the congruence $\sim$ in the preorder
$\preceq$ together with the inclusion
${\preceq} \subseteq {\preceq_i}$ imply the inclusion
${\sim} \subseteq {\preceq_i}$.
As $\sim$ is a congruence, this last inclusion yields
${\sim} \subseteq {\sim_i}$.
Thus, ${\sim} \subseteq (\bigcap_{i \; \in \; I} \sim_i)$.

2) If a congruence $\sim$ is c-residually finite and 
${\sim} = \bigcap_{i \; \in \; I} \sim_i$ then the congruences
$\sim_i$'s, being also preorders, witness that $\sim$ is
sp-residually finite. 
Conversely, applying item 1 to a congruence $\sim$, we see
that if $\sim$ is sp-residually finite then it is also c-residually
finite.
\end{proof}
The following proposition gives very simple examples.
\begin{proposition}\label{p res fin algebras}
$\langle \N; + \rangle$ (i.e. the free monoid with one
generator) is sp-residually finite whereas
$\langle \Z; + \rangle$ (i.e. the free group with one
generator) is c-residually finite but not sp-residually finite.
\end{proposition}
\begin{proof}
{\it$\langle \Z; + \rangle$ is c-residually finite.}
Every congruence which is not the identity is a modular
congruence hence has finite index.
As for the identity, it is the intersection of all modular
congruences.

{\it$\langle \Z; + \rangle$ is not sp-residually finite.}
Let's see that the usual order $\leq$ is not sp-residually
finite.
Indeed, let $\sqsubseteq$ be a stable preorder which strictly
contains $\leq$ (which is the case if $\sqsubseteq$ has
finite index and contains $\leq$).
We prove that $\sqsubseteq$ is trivial: it identifies all
integers.
Since $\sqsubseteq$ strictly contains $<$, there exists
$a \in \Z$ and $\ell \neq 0$ such that
$a \sqsubseteq a+\ell$ and $a+\ell \sqsubseteq a$.
We can suppose $\ell > 0$ since otherwise it suffices to add
$-\ell$ to these inequalities. 
Adding $-\ell+1$ to the second inequality we get
$a+1 \sqsubseteq a+1-\ell$.
Since $a+1-\ell \leq a$ and $\sqsubseteq$ contains $\leq$
we also have  $a+1-\ell \sqsubseteq a$ hence by transitivity
$a+1 \sqsubseteq a$.
By stability we get $x+1 \sqsubseteq x$ for all $x$ and by
transitivity $x > y \Rightarrow x \sqsubseteq y$ for all $x, y$.
Since $\sqsubseteq$ contains $\leq$ hence also $<$, this
shows that $\sqsubseteq$ identifies all integers.

{\it$\langle \N; + \rangle$ is sp-residually finite.}
Let $\sqsubseteq$ be a stable preorder.
If it has finite index it is trivially sp-residually finite.
So, we now suppose that $\sqsubseteq$ has infinite index.
\begin{claim}
At least one of the preorders $\sqsubseteq$ and its reverse
$\sqsupseteq$ is included in the natural order $\leq$.
In other words, at least one of the following two properties
holds:
\begin{equation}
\label{eq sqsubseteq in leq or geq infinite index}
\forall x, y \in \N \ ( x \sqsubseteq y \Rightarrow x \leq y)
\quad,\quad
\forall x, y \in \N \ ( x \sqsupseteq y \Rightarrow x \leq y).
\end{equation}
If both are true then $\sqsubseteq$ coincides with the
identity relation on $\N$.
\end{claim}
\begin{proof}[Proof of Claim]
We argue by way of contradiction.
If the two stated properties both failed, let $d, e > 0$ be
minimum such that there exist $x, y \in \N$ satisfying
$x \sqsubseteq x+d$ and $y+e \sqsubseteq y$.
Let $z = \max(x, y)$. Stability yields
$z+e \sqsubseteq z \sqsubseteq z+d$ (just add $z-x$ and
$z-y$). 

Let $f = gcd(d, e)$.
Applying Bachet-B\'ezout's theorem, 
$\alpha d - \beta e = \gamma e - \delta d = f$ for some
$\alpha, \beta, \gamma, \delta \in \N$.
Adding $d$ and $e$ again and again, we get
\begin{align*}
&z \sqsubseteq z+d 
  \sqsubseteq z+2d 
  \sqsubseteq \ldots 
  \sqsubseteq z+\alpha d = z+\beta e +f
\textit{\quad and \quad}
z \sqsupseteq z+e 
  \sqsupseteq z+2e 
  \sqsupseteq \ldots 
  \sqsupseteq z+\beta e
\\
&z \sqsubseteq z+d 
  \sqsubseteq z+2d 
  \sqsubseteq \ldots 
  \sqsubseteq z+\delta d = z+\gamma e -f
\textit{\quad and \quad}
z \sqsupseteq z+e 
  \sqsupseteq z+2e 
  \sqsupseteq \ldots 
  \sqsupseteq z+\gamma e
\end{align*}
Thus, 
$z+\beta e \sqsubseteq z \sqsubseteq z+\beta e +f$
and 
$z+\gamma e \sqsubseteq z \sqsubseteq z+\gamma e -f$.
Letting $x' = z+\beta e$ and $y' = z+\gamma e$, we get
$x' \sqsubseteq x'+f$, and (adding $f$)
$y'+f \sqsubseteq y'$.
Minimality of $d, e$ insures that $f \geq d$ and $f \geq e$.
Letting $z'=\max(x',y')$ stability yields
 $z'+f \sqsubseteq z' \sqsubseteq z'+f$.
Since $\sqsubseteq$ has infinite index its associated
congruence is the identity (cf. Lemma~\ref{l folklore cong})
so that $z' = z'+f$ hence $f = 0$, contradicting inequalities
$d, e > 0$, $f \geq d$, and $f \geq e$.
\end{proof}
Using the above Claim, we can now finish the proof.
\\\textbullet\
{\it Case $\forall x, y \in \N \ ( x \sqsubseteq y \Rightarrow
 x \leq y)$ holds.}
For $i\in\N$ let 
\begin{align}\label{eq prceqi}\tag{*}
\preceq_{i} &=\ \sqsubseteq \cup \{(x,y) \mid y\geq i\}
\end{align}
Let's check that $\preceq_{i}$ is a stable preorder.
Reflexivity is trivial. As for transitivity, if $x \preceq_{i} y$ and
$y \preceq_{i} z$ then the Case assumption insures that
$x \leq y \leq z$.
If $ z\geq i$ then (*) insures $ x\preceq_{i} z$.
If $z < i$ then also $y < i$ and (*) insures $x \sqsubseteq y$
and $y \sqsubseteq z$ hence $x \sqsubseteq z$ by
transitivity of $\sqsubseteq$.
Observe also that $\preceq_{i}$ has finite index (at most
$i+1$).

Since $\sqsubseteq \, = \bigcap_{i\in\N} \preceq_{i}$ we see
that $\sqsubseteq$ is residually finite.
\\\textbullet\
{\it Case $\forall x, y \in \N \ ( x \sqsubseteq y \Rightarrow
 x \geq y)$ holds.}
Similar, letting 
$\preceq_{i} \ =\ \sqsubseteq \cup \{(x,y) \mid x\geq i\}$.
\end{proof}
Usual residual finiteness is strictly weaker than c-residual
finiteness.
\begin{proposition}\label{p res fin not c-res fin a}
Though free groups and free monoids are residually finite,

1. Free groups on at least two generators (possibly infinitely
many) are not c-residually finite.

2. Free monoids with at least four generators (possibly
infinitely many) are not c-residually finite.
\end{proposition}
\begin{proof}
1.
G. Baumslag proved (cited in Magnus, 1969 \cite{Magnus69}
pp. 307--308) that (denoting $e$ the empty word) the
quotient of the free group $F_{2}$ with the two generators
$a$ and $b$ by the relation $a^{-1}b^{2}ab^{3} = e$ is not
residually finite.
Which means that the congruence $\approx$ of $F_{2}$
generated by the relation $a^{-1}b^{2}ab^{3} = e$ is not
c-residually finite. 
Hence $F_{2}$ is not c-residually finite.

For free groups with more generators than the sole $a$ and
$b$, simply add to the previous relation the relations
equating to $e$ all generators different from $a$ and $b$.
 
2. 
Recall that if $\theta \colon \+A \to \+B$ is a surjective
homomorphism between the algebras $\+A$ and $\+B$
and $\triplesim$ is a congruence on $\+B$ then 
$\theta^{-1}(\triplesim)$ is a congruence on $\+A$ and the
quotient algebras $\+A/\theta^{-1}(\triplesim)$ and
$\+B/\triplesim$ are isomorphic.
 
Let $\Sigma = \{ \ell_{1}, \ell_{2}, \ell_{3}, \ell_{4}\}
 \cup \Delta$ and let $\theta \colon \Sigma^* \to F_{2}$ be
the surjective semigroup morphism such that ($a$ and $b$
being the two generators of $F_{2}$)
$\theta(\ell_{1}) = a$, $\theta(\ell_{2}) = a^{-1}$, 
$\theta(\ell_{3}) = b$, $\theta(\ell_{4}) = b^{-1}$,
and $\theta(x)$ is arbitrarily chosen in $F_{2}$ for every
$x \in \Delta$.
Using item 1, we can choose a non-residually finite
congruence $\triplesim$ on $F_{2}$.
Since $F_{2}/\triplesim$ is then not residually finite, neither
is $\Sigma^*/\theta^{-1}(\triplesim)$.
Thus, the congruence $\theta^{-1}(\triplesim)$ on
$\Sigma^*$ is not c-residually finite.
\end{proof}

%
\subsection{sp/c-residually finite algebras and the
upper/lower trapeziums}
\label{ss improving for c-residually finite}
%
The above notions of c-residually finiteness and
sp-residually finiteness easily lead to the following version
of Theorem~\ref{th:LatticeGen}.
\begin{theorem}\label{thm residually finite utile}
Let $\+A = \langle A; \Xi \rangle$ be an algebra and
$f \colon A \to A$.

1. If $\+A$ is sp-residually finite and $f$ preserves all finite
index stable preorders then $f$ is stable preorder preserving 
and congruence preserving. Thus, 
\begin{equation*}
(\textit{PreordPres})^{\finind}_{\+A} 
\Longleftrightarrow 
(\textit{PreordPres})_{\+A}
\text{ and }
(\textit{CongPres})^{\finind}_{\+A} 
\Longleftrightarrow 
(\textit{CongPres})_{\+A}
\end{equation*}
and the upper and lower horizontal trapeziums in
Table~\ref{table equivalences} become series of logical
equivalences.

2. If $\+A$ is c-residually finite and $f$ preserves all finite
index congruences then $f$ is congruence preserving. Thus, 
$(\textit{CongPres})^{\finind}_{\+A} 
\Longleftrightarrow 
(\textit{CongPres})_{\+A}$
and the lower trapezium in Table~\ref{table equivalences}
becomes a series of logical equivalences.
\end{theorem}
\begin{proof}
1) Let $\preceq$ be a stable preorder. 
The hypothesis of sp-residual finiteness of $\+A$ ensures
that $\preceq$ is sp-residually finite: there exists a family of
stable preorders $(\preceq_i)_{i\; \in \; I}$ with associated
congruences having finite indexes, such that
${\preceq} = {\bigcap_{i\; \in \; I} \preceq_i}$.
Thus, $a \preceq b$ if and only if, for all $i \in I$,
$a \preceq_i b$.  
The hypothesis ensures that $f$ preserves the $\preceq_i$'s
hence  $f(a) \preceq_i f(b)$ for all $i \in I$ and therefore
$f(a) \preceq f(b)$.
Lemma~\ref{l sp implies c} leads to the stated result about
congruence preservation.
The proof of 2) is similar.
\end{proof}
As a corollary of Theorem~\ref{thm residually finite utile},
we get a new proof of the equivalence of the first and third
conditions in Theorem~\ref{thm:stable preserv and lattices}.
\begin{theorem}\label{thm:ipl1Improv}
Let $\+N$ be $\langle \N; \suc \rangle$, 
$\langle \N; + \rangle$ or 
$\langle \N; +, \times \rangle$.
The following conditions are equivalent for any function
$f\colon \N \longrightarrow \N$.
\\
\begin{equation*}
\begin{array}{ll}
(\textit{PreorderPres})_{\N}
&
\text{$f$ is $\+N$-stable-preorder preserving}
\smallskip
\\
(\Latt^{\infty}_{\freez^{*}}(\forall L))_{\+A}
& \textit{$f^{-1}(L) \in \Latt_{\+A}^{\infty}(L)$
          for every recognizable $L \subseteq A$}
\end{array}
\end{equation*}
\end{theorem}

\subsection{Residually finite algebra with a group operation}
%
When there is a group operation in a c-residually finite
algebra, Theorem \ref{th:LatticeGen} can be strengthened to 
the equicalence of all considered conditions.
\begin{theorem}\label{th:RecLatGeneralGroup} 
Let $\+A = \langle A; \Xi \rangle$ be a  c-residually finite
algebra such that $\Xi$ contains a group operation (for
instance $\langle \Z; +, \times \rangle$).  
Let $f \colon A\to A$. 
The eight conditions of Theorem~\ref{th:LatticeGen} are
equivalent.
In other words all arrows in the diagram of
Table~\ref{table equivalences} are bi-implications.
\end{theorem}
\begin{proof}
If in the horizontal upper, lower trapeziums and the vertical
right-sided trapezium the implication arrows can be
replaced by bi-implication arrows,  so also does the vertical
left-sided trapezium.
\end{proof}
%


\end{document}